\documentclass[a4paper,reqno,10pt,oneside]{amsart}
\usepackage{amsmath,geometry}
\usepackage{graphicx}
\usepackage{mathrsfs}
\usepackage{amssymb,amsmath, amsfonts, amsthm}
\usepackage{latexsym}
\usepackage{color} 
\usepackage[dvips]{epsfig}
\usepackage{graphicx}

\allowdisplaybreaks[1]

\numberwithin{equation}{section}

\renewcommand{\(}{\left(}
\renewcommand{\)}{\right)}
\renewcommand{\[}{\left[}
\renewcommand{\]}{\right]}

\newtheorem{theorem}{Theorem}[section]
\newtheorem{proposition}[theorem]{Proposition}
\newtheorem{corollary}[theorem]{Corollary}
\newtheorem{lemma}[theorem]{Lemma}
\newtheorem{remark}[theorem]{Remark}

\newcommand{\A}{{\mathcal A}}
\newcommand{\F}{{\cal F}}

\newcommand{\RN}{{\mathbb{R}^N}}

\newcommand{\N}{\mathbb{N}}

\newcommand{\beq }{\begin{equation}}
\newcommand{\eeq }{\end{equation}}

\def\F{{\mathcal F}}

\newcommand{\beqs}{\begin{equation*}}
\newcommand{\eeqs}{\end{equation*}}
\newcommand{\beqn}{\begin{eqnarray}}
\newcommand{\eeqn}{\end{eqnarray}}
\newcommand{\beqns}{\begin{eqnarray*}}
\newcommand{\eeqns}{\end{eqnarray*}}
\newcommand{\bdoc}{\begin{document}}
\newcommand{\edoc}{\end{document}}
\newcommand{\be}{\begin{enumerate}}
\newcommand{\ee}{\end{enumerate}}
\newcommand{\bdescr}{\begin{description}}
\newcommand{\edescr}{\end{description}}
\newcommand{\ba}{\begin{array}}
\newcommand{\ea}{\end{array}}
\newcommand{\intR}{\int_{\mathbb R^N}}
\newcommand{\R}{\mathbb R}
\newcommand{\Mm}{\mathcal{M}_{\lambda,\Lambda}^-}
\newcommand{\Mp}{\mathcal{M}_{\lambda,\Lambda}^+}
\newcommand{\Mpm}{\mathcal{M}_{\lambda,\Lambda}^\pm}
\newcommand{\Nm}{\tilde N_-}
\newcommand{\Np}{\tilde N_+}
\newcommand{\Npm}{\tilde N_\pm}
\newcommand{\psp}{p^*_+}
\newcommand{\psm}{p^*_-}
\newcommand{\pspm}{p^*_\pm}
\newcommand{\psc}{p^{**}_+}
\newcommand{\psce}{p_n}
\newcommand{\parallelsum}{\mathbin{\!/\mkern-5mu/\!}}
\newcommand{\e}{\varepsilon}
 \renewcommand{\(}{\left(}
\renewcommand{\)}{\right)}
\renewcommand{\[}{\left[}
\renewcommand{\]}{\right]}
\newenvironment{Proof}{\noindent{\bf Proof}}{\hfill$\Box$\\[2mm]}
\begin{document}
\title[New concentration phenomena for radial fully nonlinear equations]{New concentration phenomena for   a class of radial fully nonlinear equations}
\author{G\MakeLowercase{iulio} Galise, A\MakeLowercase{lessandro} Iacopetti, F\MakeLowercase{abiana} Leoni, F\MakeLowercase{ilomena} Pacella}

\thanks{\emph{Acknowledgements.} This research is partially supported by INDAM-GNAMPA }
\address{Dipartimento di Matematica, Sapienza Universit\`a di Roma,  P.le Aldo Moro 2, €"00185 Roma, Italy}
\email{galise@mat.uniroma1.it (G. Galise),}
\email{iacopetti@mat.uniroma1.it (A. Iacopetti),}
\email{leoni@mat.uniroma1.it (F. Leoni),}
\email{pacella@mat.uniroma1.it (F. Pacella).}

\begin{abstract}
We study radial sign-changing solutions of a class of fully nonlinear elliptic Dirichlet problems in a ball, driven by the extremal Pucci's operators and with a power nonlinear term. We first determine a new critical exponent related to the existence or nonexistence of such solutions. Then we analyze the asymptotic  behavior of the radial nodal solutions as the exponents approach the critical values, showing that new concentration phenomena occur. Finally we define a suitable weighted energy for these solutions and compute its limit value.
\end{abstract}

\maketitle
\section{Introduction}

Let $B$ be the unit ball of $\R^N$ and let $0<\lambda\leq \Lambda$.  We consider the  problem
\begin{equation}\label{eq:probMgen}
\begin{cases}
-\F(D^2 u)=|u|^{p-1}u & \hbox{in} \ B\\
\qquad  \ \ \ \ \ u=0 & \hbox{on} \ \partial B\\
\qquad  \ \ \ \ \ u(0)>0 &
\end{cases}
\end{equation}
where $p>1$, $\mathcal{F}$ is either one of the Pucci's extremal operators $\Mpm$, defined respectively as

\begin{equation*}
\begin{split}
{\mathcal M}^-_{\lambda,\Lambda}(X)&:=\inf_{ \lambda I\leq A\leq \Lambda I}{\rm tr} (AX)=\lambda\, \sum_{\mu_i>0}\mu_i+\Lambda\, \sum_{\mu_i<0}\mu_i \\
{\mathcal M}^+_{\lambda,\Lambda}(X)&:=\sup_{ \lambda I\leq A\leq \Lambda I} {\rm tr} (AX)=\Lambda\, \sum_{\mu_i>0}\mu_i+\lambda\, \sum_{\mu_i<0}\mu_i,
\end{split}
\end{equation*}
$\mu_1,\ldots,\mu_N$ being the eigenvalues of any squared symmetric matrix $X$.\\

Obviously when $\lambda=\Lambda$, \eqref{eq:probMgen} is the classical Lane-Emden problem, because  Pucci's operators reduce to a multiple of the Laplacian.

Let us immediately observe that, since ${\mathcal M}^+_{\lambda,\Lambda}(-X)=-{\mathcal M}^-_{\lambda,\Lambda}(X)$, solutions of \eqref{eq:probMgen}
for $\F=\Mm$ are  solutions of the analogous problem for the operator $\F=\Mp$, but with $u(0)<0$. Thus, it is important to fix the sign at the center of the ball.\\

The study of \eqref{eq:probMgen}, apart from being interesting in itself, is important to understand some invariance of the Pucci's operators which may be not so evident by their definition.

Indeed, though \eqref{eq:probMgen} does not have a variational structure when $\lambda<\Lambda$ (as it happens instead for the classical Lane-Emden problem) some critical exponents appear in connection with the existence of solutions. For positive solutions (which are radial by the symmetry result of \cite{DS}) they are related to the existence of radial fast decaying solutions of the analogous problem in $\R^N$ (see \cite{FQ}) and induce a concentration phenomenon for positive solutions of \eqref{eq:probMgen}, as $p$ approaches the critical values (see \cite{BGLP}). Moreover a weighted related ``energy'' was defined in \cite{BGLP} which is preserved in the limit, thought the positive solutions concentrate at the origin and converge to zero everywhere else.\\

The aim of the present paper is to study the asymptotic behavior of radial sign-changing solutions of \eqref{eq:probMgen} as the exponent $p$ approaches some critical values for their existence. We will show that new critical exponents and new concentration phenomena occur, quite different from those related to the classical Lane-Emden problem but also different from those shown in \cite{BGLP} for the positive solutions of \eqref{eq:probMgen}.

First of all we prove that a new critical exponent $\psc$ appears for the existence of radial nodal solutions to \eqref{eq:probMgen} when $\F=\Mp$, which is in between those for the existence of radial positive solutions for the two Pucci's operators (see \eqref{eq2mainteo1}). This is somehow surprising because, since the solutions of \eqref{eq:probMgen} are positive in the first nodal region, which is a ball, one would expect the critical exponents to be the same as the one for positive solutions to \eqref{eq:probMgen}. Indeed this is the case for $\F=\Mm$ and for the classical Laplacian, but not for $\F=\Mp$ (see Theorem \ref{mainteo1critexp}). 

Then we perform an accurate asymptotic analysis of radial nodal solutions of \eqref{eq:probMgen} with any number $k$ of nodal domains and show that the behavior can be different in each nodal region and may also depend on $k$ being even or odd (see Theorem \ref{mainteo2} and Theorem \ref{mainteo4}). Indeed while in some nodal domain there is blow up and concentration in others the solutions are bounded and converge to a finite limit. Moreover the asymptotic profile of the solutions $u$ of \eqref{eq:probMgen}, after suitable rescalings, can be different and, in the case of $\F=\Mp$, the fast decaying radial positive solution of \eqref{eq:domextmmpsc} in the exterior of the ball appears as limit profile of the restriction of $u$ to some nodal regions (see Proposition \ref{prop:unifupperbound}).

This is a completely new phenomenon, to our knowledge, different from what happens for the classical Lane-Emden problem (see \cite{DMIP} and the references therein) and even from what happens in the case of the classical Brezis-Nirenberg problem in low-dimensions which also presents some peculiar asymptotic behavior (see \cite{ABP}, \cite{AP}, \cite{Iac}, \cite{IacPac2}, \cite{IacVair}, \cite{IacVair2}).

Finally, all this reflects into the computation of the limit of some weighted energies which can be defined for solutions of \eqref{eq:probMgen}, according to what done in \cite{BGLP}, even if \eqref{eq:probMgen} does not have a variational structure.

We will show that the weighted energy of the positive fast decaying solutions, both in $\R^N$ and in $\R^N\setminus B$ will contribute to the limit of the total energy of $u$ in some of the nodal regions where blow up and concentration occur.

To state precisely our results let us start by recalling what is known for positive solutions to \eqref{eq:probMgen}.

In the paper \cite{FQ} Felmer and Quaas proved that there exist two critical exponents $\psm$, $\psp$ such that positive radial classical solutions to \eqref{eq:probMgen} exist if and only if  $p<\psm$ for $\F=\Mm$ or $p<\psp$ when $\F=\Mp$. We observe that the values of these critical exponents are not explicitly known but they satisfy the following inequalities:
\begin{equation}\label{criticalexlimitations}
\begin{split}
\frac{\Nm+2}{\Nm-2}&<\,\psm<\frac{N+2}{N-2},\\
\max\left\{\frac{\Np}{\Np-2}, \frac{N+2}{N-2}\right\}&<\,\psp<\frac{\Np+2}{\Np-2},
\end{split}
\end{equation}
where the dimension-like parameters $\Npm$ are defined, respectively by $\tilde N_-:=\frac{\Lambda}{\lambda}(N-1) +1$, $\tilde N_+:=\frac{\lambda}{\Lambda}(N-1) +1$.

We point out that in the special case $\lambda=\Lambda$, when $\Mm=\Mp=\lambda\Delta$, where $\Delta$ is the standard Laplacian, all the above inequalities become equalities. In particular $\psm$, $\psp$ reduce to the usual Sobolev critical exponent $\frac{N+2}{N-2}$.
 
In this paper we will always assume that $\lambda<\Lambda$ and $\Npm>2$.\\ 

As far as the existence  of radial sign-changing solutions is concerned, let us mention that in \cite{GLP} a sufficient condition on the exponent  $p$ is provided for general  radially symmetric nonlinear operators which, in the particular case of the Pucci's operators, reads as $p\leq\frac{\Nm}{\Nm-2}$.\\ 

The first result of the present paper shows that such a bound on the exponent $p$ is not optimal. Indeed we get:

\begin{theorem}\label{mainteo1critexp}
We have:
\begin{enumerate}
\item[i)] if $\F=\Mm$, then radial sign-changing solutions of \eqref{eq:probMgen} with any number of nodal domains exist if and only if
\beq\label{eq1mainteo1} 
p<\psm\,;
\eeq
\item[ii)] if $\F=\Mp$, then there exists a new critical exponent $\psc$ satisfying
\beq\label{eq2mainteo1} 
\psm<\psc<\psp\,,
\eeq
such that no radial sign-changing solutions to \eqref{eq:probMgen} exist for $p\geq \psc$, while radial sign-changing solutions to \eqref{eq:probMgen} with any number of nodal domains exist at least for a sequence of exponents $p_n\nearrow\psc$.
\end{enumerate}
\end{theorem}
The above result will be proved in Section 3. Let us observe that while it is easy to obtain i), using Theorem 3.1 of \cite{GIL}, the proof of ii) is quite involved and requires several steps.\\

Once we have these critical exponents we proceed studying the asymptotic behavior of the nodal solutions of \eqref{eq:probMgen} as $p$ approaches them to determine also their limit profile. As announced before, we will see that new concentration phenomena occur.\\

We first start by analyzing the case when $\F=\Mm$. Let $p_\e:=\psm-\e$, where $\e>0$ is a small parameter and let us consider the problem
\begin{equation}\label{eq:probM-}
\begin{cases}
-\Mm(D^2 u)=|u|^{p_\e-1}u & \hbox{in} \ B\\
\qquad  \ \ \ \ \ u=0 & \hbox{on} \ \partial B\\
\qquad  \ \ \ \ \ u(0)>0 &
\end{cases}
\end{equation}
Let $u_\e$ be a radial sign-changing solution of \eqref{eq:probM-} with $k\geq 2$ nodal regions. We denote by $r_1=r_1(\e)<\ldots<r_{k-1}=r_{k-1}(\e)$ the nodal radii of $u_\e$ and by $s_i=s_{i}(\e)$  the unique maximum points of $|u_\e|$ in the $(i+1)$-th nodal region, for $i=0,\ldots,k-1$. We have $$0=s_0<r_1<s_1<\ldots<r_{k-1}<s_{k-1}<1,$$ 
and we set $M_i:=|u_\e(s_i)|$, for $i=0,1,\ldots,k-1$.

\begin{theorem}\label{mainteo2}
Up to a subsequence, as $\e \to 0^+$, we have that $M_0\to +\infty$, $r_1\to 0$, $s_1\to 0$, $M_i \to \bar M_i \in (0,+\infty)$, for $i=1,\ldots,k-1$, and 
$u_\e \to \bar u$ in $C^2_{loc}(\overline B\setminus\{0\})$, where $\bar u$ is a radial sign-changing solution of
\begin{equation}\label{eq:probM-crit2}
\begin{cases}
-\Mm(D^2 u)=|u|^{p^*_- -1}u & \hbox{in} \ B\\
\qquad  \ \ \ \ \ u=0 & \hbox{on} \ \partial B\\
\qquad  \ \ \ \ \ u(0)<0 &
\end{cases}
\end{equation}
with $(k-1)$ nodal regions, if $k\geq 3$, while $\bar u$ is the unique negative solution of \eqref{eq:probM-crit2} if $k=2$.

Moreover if $k\geq 3$ we have $r_i \to \bar r_i$, $s_i \to \bar s_i$, for  $i=2,\ldots,k-1$, for some numbers $\bar r_i, \bar s_i$, such that $0<\bar r_2<\bar s_2<\ldots<\bar r_{k-1}<\bar s_{k-1}<1$.
\end{theorem}

Note that \eqref{eq:probM-crit2} does not admit a positive solution, by \eqref{criticalexlimitations}, but it has a (unique) negative solution as well as sign-changing solutions by \eqref{eq1mainteo1}, since $u(0)<0$ so that the relevant exponents for the existence of solutions to \eqref{eq:probM-crit2} are those for the corresponding equations involving the operator $\Mp$, but requiring the positivity at the origin.\\

Even if the problems that we are considering do not have a variational structure, we can introduce, in the spirit of \cite{BGLP},  a weighted energy $E^T_p(u)$ defined for radial sign-changing functions $u$ which change concavity only once in each nodal region, and where $p>1$ is a fixed exponent (we refer to Sect. 10 for the definition). In particular, if $u_\e$ is as in the statement of Theorem \ref{mainteo2} we are interested in determining the limit energy $E_{p_\e}^T(u_\e)$ as $\e\to 0^+$. To this end, denoting by $U_-$ the unique (up to scaling) positive radial fast decaying solution of
$$-\Mm(D^2u)=u^{\psm} \ \ \text{in} \ \R^N,$$
and setting
\beq\label{limitenergysigmaU-}
\Sigma^*_-:=E^*(U_-),
\eeq
where $E^*(U_-)$ is the (finite) energy of $U_-$ in $\RN$ with $p=\psm$ (see \eqref{weightedenergyRN}) we have the following.
\begin{theorem}\label{mainteo3}
Let $u_\e$ be as in Theorem \ref{mainteo2}. It holds
\begin{equation}\label{eq:limenergy}
\lim_{\e \to 0^+} E_{p_\e}^T(u_\e) = \Sigma^*_- + E_{\psm}^T(\bar u),
\end{equation}
where $\Sigma^*_-$ is defined by \eqref{limitenergysigmaU-} and $E_{\psm}^T(\bar u)$ is the total energy of the limit function $\bar u$ (given by Theorem \ref{mainteo2}), i.e.
\begin{equation}\label{totalenergylimit}
 E_{\psm}^T(\bar u) = \sum_{j=1}^{k-1}E_{\psm, {\Omega}^j}(\bar u^j), \ \ j=1,\ldots,k-1,
\end{equation}
where $\bar u^j$ is the restriction of $\bar u$ to its j-th nodal region ${\Omega}^j$, $j=1,\ldots, k-1$ and $E_{\psm, {\Omega}^j}(\bar u^j)$ is its energy as defined in \eqref{restrweightedenergy}.
\end{theorem}

When $\F=\Mp$ the picture is quite different. Setting $p_n:=\psc-\e_n$, where $\e_n>0$ is a sequence converging to zero as $n\to +\infty$, we consider a radial sign-changing solution $u_n$ of the problem 
\begin{equation}\label{eq:probM+}
\begin{cases}
-\Mp(D^2 u)=|u|^{\psce-1}u & \hbox{in} \ B\\
\qquad  \ \ \ \ \ u=0 & \hbox{on} \ \partial B\\
\qquad  \ \ \ \ \ u(0)>0 &
\end{cases}
\end{equation}
As before, for $k\geq 2$ and $i=1,\ldots,k-1$, we denote by $r_i=r_i(n)$, the nodal radii of $u_n$, by $s_i=s_i(n)$ the unique maximum point in the $(i+1)$-th nodal region and define $M_i=|u_n(s_i)|$.
\begin{theorem}\label{mainteo4}
Up to a subsequence, as $n\to +\infty$, we have:
\begin{itemize}
\item[i)] if $k$ is even then $M_0\to +\infty$, $M_i\to +\infty$, $r_i\to 0$, $s_i\to 0$ for all $i=1,\ldots,k-1$, and $u_n \to 0$ in $C^2_{loc}(\overline{B}\setminus\{0\})$. Moreover for $j=0,\ldots,\frac{k-2}{2}$ there exist positive constants $c_j$ such that  $$\frac{M_{2j}}{M_{2j+1}}\to c_j;$$ if $k\geq4$ we also have that $$\frac{M_{2j+1}}{M_{2j+2}}\to+\infty,$$ for $j=0,\dots,\frac{k-4}{2}$, as $n\to +\infty$; 
\item[ii)] if $k$ is odd then $M_0\to +\infty$, $M_i\to +\infty$, $r_i\to 0$, $s_i\to 0$ for all $i=1,\ldots,k-2$, $r_{k-1} \to 0$, $s_{k-1} \to 0$, $M_{k-1} \to \bar M$, for some $\bar M>0$ and $u_n \to \bar v$ in $C^2_{loc}(\overline{B}\setminus\{0\})$, where $\bar v$ is the unique positive solution of
\beq\label{probllimupodd}
\begin{cases}
-\Mp(D^2 u)=u^{\psc} & \hbox{in} \ B\\
\quad \quad \quad \quad\quad \  u= 0& \hbox{on} \ \partial B.
\end{cases}
\eeq

Moreover for $j=0,\ldots,\frac{k-3}{2}$ there exist positive constants $c_j$ such that, as $n\to +\infty$, 
\begin{equation*}
\frac{M_{2j}}{M_{2j+1}} \to c_j,\ \ \frac{M_{2j+1}}{M_{2j+2}} \to+\infty.
\end{equation*} 
\end{itemize}
\end{theorem}

\bigskip

To determine the limit energy of $u_n$ we denote by $W_-$ the only positive radial fast decaying solution of
\begin{equation}\label{eq:domextmmpsc}
\begin{cases}
-\Mm(D^2 u)=u^{\psc} & \hbox{in} \ \R^N\setminus \overline B\\
\qquad  \ \ \ \ \ u=0 & \hbox{on} \ \partial B
\end{cases}
\end{equation}
which exists by the results of \cite{GIL} because $\psm<\psc$.
Then, setting 
\beq\label{limitenergysigmaW-}
\Sigma^{**}_+:=E^{**}(W_-),
\eeq
where $E^{**}(W_-)$ is the (finite) energy in $\R^N\setminus B$ of $W_-$ (see \eqref{eq:energyW-}), we have the following.

\begin{theorem}\label{mainteo5}
Let $u_n$ be as in the statement of Theorem \ref{mainteo4}. We have
\begin{equation}\label{eq:mainteo5}
\lim_{n\to +\infty} E_{\psce}^T(u_n) = \begin{cases}
 \frac{k}{2} E_{\psc, B}(\bar v) + \frac{k}{2}\Sigma_+^{**} & \hbox{if $k$ is even},\\[6pt]
  \frac{k+1}{2} E_{\psc, B}(\bar v) +  \frac{k-1}{2}\Sigma_+^{**}  & \hbox{if $k$ is odd},
  \end{cases}
\end{equation}
where $\Sigma_+^{**}$ is defined by \eqref{limitenergysigmaW-}, and $E_{\psc, B}(\bar v)$ is the energy of the only radial positive solution $\bar v$ to \eqref{probllimupodd} (see \eqref{weightedenergy}).
\end{theorem}

The proofs of the above results are quite involved and combine several methods: blow up techniques and study of some limit problems, phase plane analysis for the corresponding ODE's and estimates on related pointwise energies.\\

The outline of the paper is the following. In Section 2 we recall some preliminary results on positive solutions. In Section 3 we prove Theorem \ref{mainteo1critexp}. In Section 4 and Section 5 we consider the case of solutions to \eqref{eq:probMgen} for $\F=\Mm$ with two or three nodal regions. This allows to study the case of any number $k$ of nodal domains by induction in Section 6, proving so Theorem \ref{mainteo2}. In Section 7 we study problem \eqref{eq:probMgen} for $\F=\Mp$ and solutions with two nodal regions, while in Section 8 we consider the case of three nodal domains. The proof of Theorem \ref{mainteo4} is then presented in Section 9, again by an induction argument.

Finally in Section 10 we study the total energy associated to the nodal solutions of \eqref{eq:probMgen} and prove Theorem \ref{mainteo3} and Theorem \ref{mainteo5}.

\bigskip

\section{Preliminary results on positive radial solutions}
We begin this section by recalling the known results about the asymptotic analysis of positive radial solutions to
\begin{equation}\label{eq:probMpos}
\begin{cases}
-\F(D^2 u)=u^{p} & \hbox{in} \ B\\
\qquad  \ \ \ \ \ u=0 & \hbox{on} \ \partial B
\end{cases}
\end{equation}
as $p$ approaches the critical exponent for which such solutions exist (for the proofs we refer to \cite{BGLP}). We first introduce some notation: let $\e>0$ be a small parameter and set $$p_\e:=\begin{cases}\psp- \e, & \hbox{if} \ \F=\Mp,\\ \psm- \e, & \hbox{if} \ \F=\Mm.\end{cases}$$ 
We denote by $v_{p_\e, \pm}$ the unique positive solution of \eqref{eq:probMpos}. Namely $v_{p_\e,+}$ is the only positive solution to \eqref{eq:probMpos} if $\F=\Mp$, and $v_{p_\e,-}$ is the only positive  solution to \eqref{eq:probMpos} if $\F=\Mm$. Accordingly, we denote by $r_{0,\pm}=r_{0,\pm}(\e) \in (0,1)$ the only radius such that $v_{p_\e,\pm}^{\prime\prime} (r) < 0$ for $r \in [0, r_{0,\pm})$ and $v_{p_\e,\pm}^{\prime\prime} (r)>0$ for $r \in (r_{0,\pm},1)$. Moreover, let $U_{\pm}$ be the unique positive radial solution of
\beq\label{criticalproblemMPMMRN}
-\Mpm(D^2 u)=u^{p^*_\pm} \ \ \ \hbox{in} \ \RN
\eeq
such that $U_\pm(0)=1$, and denote by $R_{0,\pm}$ the unique radius such that $U_\pm^{\prime\prime} (r) < 0$ for $r \in [0, R_{0,\pm})$, $U_\pm^{\prime\prime} (r) > 0$ for $r \in (R_0^\pm,+\infty)$. We refer to the solutions of \eqref{criticalproblemMPMMRN}, and in particular to $U_\pm$, as the fast decaying solutions, since for all $p\geq p^*_\pm$ and among all radial positive solutions of
$$
-\Mpm(D^2 u)=u^{p} \ \ \ \hbox{in} \ \RN\, ,
$$
 one has
 $$
 \lim_{r\to +\infty}r^\frac{2}{p-1}u(r)=0\quad \Longleftrightarrow \quad p=p^*_\pm\, .
 $$
\begin{proposition}\label{Prop:asanpos}
Let $v_{p_\e,\pm}$ be the unique positive solution to \eqref{eq:probMpos}. Then:
\begin{itemize}
\item[(i)] $\displaystyle \lim_{\e\to 0^+}\left\|v_{p_\e,\pm}\right\|_\infty = \lim_{\e\to 0^+} v_{p_\e,\pm}(0)= +\infty$;\\
\item[(ii)] $v_{p_\e,\pm} \to 0$ in $C^2_{loc}(\overline B \setminus \{0\})$, as $\e\to 0^+$;
\item[(iii)]  $ \displaystyle \lim_{\e\to 0^+}  [r_{0,\pm}(\e)]^{\frac{2}{p_\e - 1}}\left\|v_{p_\e,\pm}\right\|_\infty = (R_{0,\pm})^{\frac{2}{p^*_\pm - 1}}$;
\item[(iv)] $\displaystyle \lim_{\e\to 0^+}  \frac{v_{p_\e,\pm}(r_{0,\pm}(\e))}{\left\|v_{p_\e,\pm}\right\|_\infty} = U_\pm(R_{0,\pm})$;
\item[(v)] $\displaystyle \lim_{\e\to 0^+}  [r_{0,\pm}(\e)]^{\frac{2}{p_\e - 1}}v_{p_\e,\pm}(r_{0,\pm}(\e))= (R_{0,\pm})^{\frac{2}{p^*_\pm - 1}}U_\pm(R_{0,\pm})$;
\item[(vi)] $\displaystyle\lim_{\e \to 0^+} \left\|v_{p_\e,\pm}\right\|_\infty^{\frac{p_\e(\tilde N_\pm-2) - \tilde N_\pm}{2}} (v_{p_\e,\pm})^\prime(1)=-C_\pm$, where  $C_\pm$ is a positive constant depending only on $N,\lambda,\Lambda$.
\end{itemize}
\end{proposition}

Next we recall some useful results about the qualitative properties of the solutions of a suitable class of initial value problems. To do this we need some preliminaries.

If $u$ is a smooth radially symmetric function, we easily check that the Hessian of $u$ is given by
\beq\label{eq:hessianradial}
D^2u(x)=\frac{u^\prime(|x|)}{|x|} \mathbb{I}_N + \left(u^{\prime\prime}(|x|)-\frac{u^\prime(|x|)}{|x|}\right)\frac{x}{|x|}\otimes \frac{x}{|x|},
\eeq
where $\mathbb{I}_N$ is the identity matrix of order $N$ and $x\otimes x$ is the matrix defined by $(x\otimes x)_{ij}=x_ix_j$, for any $i,j\in\{1,\ldots,N\}$. In particular, since the eigenvalues of the matrix appearing in the right-hand side of \eqref{eq:hessianradial} are $u^{\prime\prime}(|x|)$, which is simple, and $\frac{u^\prime(|x|)}{|x|}$, which has multiplicity $(N-1)$, we infer that if $u$ is a positive radial solution of $-\F(D^2 u)=u^{p}$, then setting $r=|x|$ there are only three possibilities: \\

 \textbf{Case 1}: $u^\prime(r)\geq 0$ and $u^{\prime\prime}(r)\leq 0$, so that $u=u(r)$ satisfies
 \begin{equation}\label{case1}
 \begin{cases}
 -\Lambda u^{\prime\prime}(r) - \lambda (N-1)\frac{u^\prime(r)}{r}=u^p(r) & \text{if $\F=\Mm$},\\[6pt]
  -\lambda u^{\prime\prime}(r) - \Lambda (N-1)\frac{u^\prime(r)}{r}=u^p(r) & \text{if $\F=\Mp$}.
 \end{cases}
 \end{equation}
 \textbf{Case 2}: $u^\prime(r)\leq 0$ and $u^{\prime\prime}(r)\leq 0$, so that $u=u(r)$ satisfies
 \begin{equation}\label{case2}
 \begin{cases}
 -\Lambda \left(u^{\prime\prime}(r)+(N-1)\frac{u^\prime(r)}{r}\right)=u^p(r) & \text{if $\F=\Mm$},\\[6pt]
  -\lambda \left(u^{\prime\prime}(r) + (N-1)\frac{u^\prime(r)}{r}\right)=u^p(r) & \text{if $\F=\Mp$}.
 \end{cases}
 \end{equation}
 \textbf{Case 3}: $u^\prime(r)\leq 0$ and $u^{\prime\prime}(r)\geq 0$, so that $u=u(r)$ satisfies
 \begin{equation}\label{case3}
 \begin{cases}
 -\lambda u^{\prime\prime}(r) - \Lambda (N-1)\frac{u^\prime(r)}{r}=u^p(r) & \text{if $\F=\Mm$},\\[6pt]
  -\Lambda u^{\prime\prime}(r) - \lambda (N-1)\frac{u^\prime(r)}{r}=u^p(r) & \text{if $\F=\Mp$}.
 \end{cases}
 \end{equation}
We stress that the case $u^\prime(r)\geq 0$ and $u^{\prime\prime}(r)\geq 0$ cannot occur because $u>0$ satisfies $-\F(D^2 u)=u^{p}$.

Now, let $\alpha>0$, $p>1$ and consider the following initial value problem
\begin{equation}\label{eq:initvalprob-}
\begin{cases}
u^{\prime\prime}(r)=M_-\left(-\frac{\Lambda(N-1)}{r}K_-(u^\prime(r)) - u^{p}(r)\right) & \hbox{for} \ r>1\\
\qquad  \ \ \ \ \ u(r)>0 & \hbox{for} \ r>1\\
\qquad  \ \ \ \ \ u(1)=0, \ \ u^\prime(1)=\alpha& 
\end{cases}
\end{equation}
with
$$M_-(\xi):=\begin{cases} \xi/\lambda & \hbox{if} \ \xi\geq 0\\ \xi/\Lambda & \hbox{if} \ \xi< 0\end{cases}, \ \ \ \ \ K_-(\xi):=\begin{cases} \frac{\lambda}{\Lambda}\xi & \hbox{if} \ \xi\geq 0\\ \xi& \hbox{if} \ \xi< 0.\end{cases}$$

Problem \eqref{eq:initvalprob-} has a unique solution $u_\alpha=u(\alpha,p,r)$, defined and positive on a maximal interval $[1, \rho_\alpha)$, for some $1<\rho_\alpha\leq +\infty$. In \cite{GLP} it has been proved that there exists $\tau_\alpha \in (1,\rho_\alpha)$ such that $u_\alpha^\prime(r)>0$ for $r \in (1,\tau_\alpha)$, $u_\alpha^\prime(r)<0$ for $r \in (\tau_\alpha, \rho_\alpha)$. Moreover, there exists $\sigma_\alpha \in (\tau_\alpha,\rho_\alpha)$ such that $u_\alpha^{\prime\prime}<0$ in $(1,\sigma_\alpha)$ and $u_\alpha^{\prime\prime}>0$ in $(\sigma_\alpha,\rho_\alpha)$ (see \cite{FQ,GIL}).

Concerning the asymptotic properties with respect to the parameter $\alpha$, we recall that $\rho_\alpha \to 1$ as $\alpha\to +\infty$, while $\rho_\alpha\to +\infty$ as $\alpha\to 0$ (see \cite[Proposition 3.2 and Lemma 3.1]{GLP}). In particular $\rho_\alpha<+\infty$ for all sufficiently large $\alpha>0$ and thus we can define the critical slope

\beq\label{defalphastar}
\alpha_-^*=\alpha_-^*(p) := \inf \left\{\alpha>0; \ \rho_\alpha<+\infty \right\}.
\eeq
We point out that if $\rho_\alpha<+\infty$ then $u_\alpha(\rho_\alpha)=0$ and $u(x):=u_\alpha(|x|)$ is a positive radial solution of
\begin{equation*}
\begin{cases}
-\Mm(D^2 u)=u^{p} & \hbox{in} \ A_{1,\rho_\alpha}\\
\qquad  \ \ \ \ \ u=0 & \hbox{on} \ \partial A_{1,\rho_\alpha}
\end{cases}
\end{equation*}
where $A_{1,\rho_\alpha}:=\{x \in \RN; \ 1<|x|<\rho_\alpha\}$ is the annulus of radii $1$, $\rho_\alpha$, centered at the origin. If $\rho_\alpha=+\infty$ then $u(x):=u_\alpha(|x|)$ is a positive radial solution of 
\begin{equation}\label{eq:probMpm}
\begin{cases}
-\Mm(D^2 u)=u^{p} & \hbox{in} \ \R^N\setminus \overline B\\
\qquad  \ \ \ \ \ u=0 & \hbox{on} \ \partial B
\end{cases}
\end{equation}

In \cite{GIL} it has been proved that \eqref{eq:probMpm} has  positive radial solutions if and only if $p>\psm$ (see \cite[Theorem 1.1]{GIL}). 
 More precisely, we have the following (see \cite[Sect. 6 and Theorem 6.2]{GIL}):
 \begin{theorem}\label{fastMm} 
Let $u_\alpha$ denote  the maximal positive solution of \eqref{eq:initvalprob-}. One has $\alpha_-^*(p)>0$ if and only if  $p>p^*_-$  and, for such $p$,
 \begin{itemize}
 \item[(i)] for any $\alpha>\alpha_-^*(p)$ it holds that $\rho_\alpha<+\infty$;
  \item[(ii)] $\rho_{\alpha_-^*}=+\infty$ and $u_{\alpha_-^*}$ is a fast decaying solution of \eqref{eq:probMpm};
 \item[(iii)] if $p^*_-<p\leq \frac{N+2}{N-2}$, then for any $\alpha<\alpha_-^*(p)$, $u_\alpha$ is either a pseudo-slow or a slow decaying solution;
  \item[(iv)] if $p> \frac{N+2}{N-2}$, then for any $\alpha<\alpha_-^*(p)$, $u_\alpha$ is a  slow decaying solution.
  \end{itemize}
\end{theorem}

Analogous results hold for $\F=\Mp$, where one considers the initial value problem

\begin{equation}\label{eq:initvalprob+}
\begin{cases}
u^{\prime\prime}(r)=M_+\left(-\frac{\lambda(N-1)}{r}K_+(u^\prime(r)) - u^{p}(r)\right) & \hbox{for} \ r>1\\
\qquad  \ \ \ \ \ u(r)>0 & \hbox{for} \ r>1\\
\qquad  \ \ \ \ \ u(1)=0, \ \ u^\prime(1)=\alpha & 
\end{cases}
\end{equation}
where 
$$M_+(\xi):=\begin{cases} \xi/\Lambda & \hbox{if} \ \xi\geq 0\\ s/\lambda & \hbox{if} \ \xi< 0\end{cases}, \ \ \ \ \ K_+(\xi):=\begin{cases} \frac{\Lambda}{\lambda}\xi & \hbox{if} \ \xi\geq 0\\ \xi& \hbox{if} \ \xi< 0\end{cases}.$$
We refer to \cite{GIL} for the precise statements.\\

We conclude this section by proving a crucial property of the map $p\mapsto\alpha_-^*(p)$.

\begin{proposition}\label{lem:asalphastar}
The map $p\mapsto\alpha_-^*(p)$ is continuous in $(1,+\infty)$.
\end{proposition}
\begin{proof}
By Theorem \ref{fastMm}, one has $\alpha_-^*(p)\equiv 0$ for $p\in (1,p^*_-]$.

 Let us first prove that $\alpha^*_-(p)\to 0$ for $p\searrow p^*_-$. 
By contradiction, assume  that there exist $\alpha_0>0$ and a sequence $p_n \searrow \psm$ such that $\alpha_-^*(p_n)>\alpha_0$ for all $n \in \N$.  This means that, for all $n \in \N$, the initial value problem
\beq\label{intialvalueprobk}
\begin{cases}
u^{\prime\prime}(r)=M_-\left(-\frac{\Lambda(N-1)}{r}K_-(u^\prime(r)) - u^{p_n}(r)\right)& \hbox{for} \ r>1\\
\qquad  \  u(1)=0, \ u^\prime(1)=\alpha_0 &
\end{cases}
\eeq
has a solution $u_{n}$ defined and positive in the whole interval $(1,+\infty)$. Let us denote by $s_n \in (1,+\infty)$ the unique maximum point of $u_n$ and set $m_n:=u_n(s_n)$. By \eqref{case1},  the energy-like functionals
$$
H_{\Lambda,n}(r):=\frac{(u_n^\prime(r))^2}{2}+\frac{(u_n(r))^{p_n+1}}{\Lambda(p_n+1)}
$$
are nonincreasing in $[1,s_n]$. Hence, we deduce
$$
m_n^{p_n+1}\leq\frac{\Lambda(p_n+1)}{2}\alpha_0^2.
$$
Then, $(m_n)_n$ is bounded and, from \eqref{intialvalueprobk}, we infer that $u_n \to \bar u$ in $C^2_{loc}([1,+\infty))$, as $n \to +\infty$, where $\bar u$ is a solution to
$$
\begin{cases}
u^{\prime\prime}(r)=M_-\left(-\frac{\Lambda(N-1)}{r}K_-(u^\prime(r)) - u^{\psm}(r)\right)& \hbox{for} \ r>1\\
\qquad  \   u\geq0 & \hbox{for} \ r>1\\
\qquad  \  u(1)=0, \ u^\prime(1)=\alpha_0.&
\end{cases}
$$
Such function $\bar u$ cannot be identically zero in view of the initial condition $\bar u'(1)=\alpha_0>0$. Hence, $\bar u>0$ in $(1,+\infty)$ and $u(x)=\bar u(|x|)$ is a positive radial solution of \eqref{eq:probMpm} with $p=\psm$,  contradicting Theorem \ref{fastMm}.

Next, let us show that $\alpha^*_-$ is continuous in $(p^*_-,+\infty)$. For any fixed $p_0>p^*_-$, let us consider $\alpha>\alpha^*_-(p_0)$. Then,  denoting by $u=u_{\alpha, p_0}$ the unique maximal solution of the initial value problem
\beq\label{ivp}
\begin{cases}
u^{\prime\prime}(r)=M_-\left(-\frac{\Lambda(N-1)}{r}K_-(u^\prime(r)) - |u|^{p_0-1}u (r)\right)& \hbox{for} \ r>1\\
\qquad  \  u(1)=0, \ u^\prime(1)=\alpha\, , &
\end{cases}
\eeq
there exists a  $\rho_\alpha>1$ such that  $u_{\alpha,p_0}(r)>0$ in $(1,\rho_\alpha)$, $u_{\alpha,p_0}(\rho_\alpha)=0$ and $u_{\alpha,p_0}<0$ in a right neighborhood of $\rho_\alpha$. By continuous dependence on the data, for $p\to p_0$ the corresponding maximal solution $u_{\alpha,p}$ is converging to $u_{\alpha,p_0}$ in $C^2_{loc}([1,+\infty))$. Hence, $u_{\alpha,p}$ has a first zero close to $\rho_\alpha$ for $p$ close to $p_0$, meaning that $\alpha^*_-(p)\leq \alpha$. By the arbitrary choice of $\alpha>\alpha^*_-(p_0)$, we deduce that $\limsup_{p\to p_0}\alpha^*_-(p)\leq \alpha^*_-(p_0)$. 

Conversely, let us now consider $0<\alpha<\alpha^*_-(p_0)$. Then, the maximal solution $u_{\alpha,p_0}$ of problem \eqref{ivp} is positive in $(1,+\infty)$.  Let us prove that, for $p$ close to $p_0$, one has $\alpha^*_-(p)\geq \alpha$. Arguing by contradiction, let us assume  that, for a sequence $p_n\to p_0$, the corresponding maximal solutions $u_{\alpha, p_n}$ satisfy $u_{\alpha,p_n}(r)>0$ in $(1,\rho_n)$ and $u_{\alpha,p_n}(\rho_n)=0$ for some  $\rho_n>1$. Again by continuous dependence on the data, one has that  $u_{\alpha,p_n}\to u_{\alpha,p_0}$ in $C^2_{loc}([1,+\infty))$, so that $\rho_n\to+ \infty$. Moreover, for each $n$  there exists $t_n\in (1,\rho_n)$ such that $u_{\alpha,p_n}''(r)<0$ for $r\in [1,t_n)$ and $u_{\alpha,p_n}''(r)>0$ for $r\in (t_n,\rho_n]$, and the sequence $(t_n)_n$ is bounded from above and from below away from 1, since otherwise the function $u_{\alpha, p_0}$ would be globally either concave or convex in $(1,+\infty)$. Thus, possibly considering a subsequence, there exists $t_0>1$ such that $t_n\to t_0$, with $u_{\alpha, p_0}''(r)<0$ for $r\in [1,t_0)$ and $u_{\alpha, p_0}''(r)>0$ for $r>t_0$. Now, we claim that there exist positive constants $C, K>0$ independent of $n$ such that
\beq\label{eq:uniformupperboundalpha}
u_{\alpha,p_n}(r)\leq  \frac{C}{\left(r^2- t_n^2 + K\right)^{\frac{\Nm-2}{2}}} \qquad \hbox{for} \ r \in [t_n,\rho_n]\, .
\eeq
Indeed, we observe that in the interval $[t_n,\rho_n]$, by \eqref{case3},  the function  $v_n=u_{\alpha,p_n}$ satisfies
\beq\label{ODEunifuppbound}
v_n^{\prime\prime}+\frac{\Nm-1}{r}v_n^\prime + \frac{v_n^{\psce}}{\lambda} =0\, .
\eeq
Then, considering the energy-like functional $H_n: [t_n, \rho_n]\to \R$ defined by
$$H_n(r):= r^{\Nm}\left[\frac{(v_n^\prime(r))^2}{2}+\frac{\Nm-2}{2\lambda \Nm} v_n^{\psce+1}\right]+ \frac{\Nm-2}{2}r^{\Nm-1}v_n(r) v_n^\prime(r),$$
and exploiting \eqref{ODEunifuppbound}, we see that
$$H_n^\prime(r)=\frac{\psce(\Nm-2)-(\Nm+2)}{2\lambda \Nm} r^{\Nm} v_n^{\psce}(r) v_n^\prime(r) <0,$$
where we use the fact that $p_n\to p_0$ and $p_0>p^*_->\frac{\tilde N_-+2}{\tilde N_--2}$.
Hence, $H_n$ is decreasing and, in particular, we get that
$$H_n(r)\geq H_n(\rho_n)=\rho_n^{\Nm}\frac{ (v_n^\prime(\rho_n))^2}{2}>0, \ \ \hbox{for any $r \in [t_n, \rho_n]$}.$$
Now, let us consider the auxiliary functional $J_n: [t_n, \rho_n]\to \R$ defined by
$$J_n(r):=v_n(r)^{-\frac{\Nm}{\Nm-2}} \frac{v_n^\prime(r)}{r}.$$
Then, exploiting again \eqref{ODEunifuppbound} and the definition of $H_n$, we easily check that
$$J_n^\prime(r)= -\frac{2\Nm}{\Nm-2} v_n(r)^{-\frac{2(\Nm-1)}{\Nm-2}} r^{-(\Nm+1)}H_n(r)<0.$$
Therefore, $J_n$ is monotone decreasing and thus $J_n(t_n) \geq J_n(r)$ for any $r\in [t_n, \rho_n]$. With the help of \eqref{ODEunifuppbound}, this can be rewritten as
$$
- \frac{v_n( t_n)^{\psce-\frac{\Nm}{\Nm-2}}}{\lambda(\Nm-1)} \geq v_n(r)^{-\frac{\Nm}{\Nm-2}} \frac{v_n^\prime(r)}{r}.
$$
Since $v_n(t_n)=u_{\alpha, p_n}(t_n)\to u_{\alpha,p_0}(t_0)>0$, 
from  the above inequality we deduce that there exists $C_1>0$ independent of $n$ such that, for any $r\in [t_n, \rho_n]$,
$$ v_n(r)^{-\frac{\Nm}{\Nm-2}} {v_n^\prime(r)} \leq -C_1 r.$$
Integrating between $t_n$ and $r$ we obtain
$$v_n(r)^{-\frac{2}{\Nm-2}} - v_n(t_n)^{-\frac{2}{\Nm-2}}\geq  \frac{C_1}{\Nm-2} (r^2-t_n^2). $$
Taking into account  that $v_n(t_n)\to u_{\alpha,p_0}(t_0)>0$ as before,  we infer that
$$ v_n(r)^{-\frac{2}{\Nm-2}} \geq  \frac{C_1}{\Nm-2} (r^2-t_n^2)+ K_1,$$ 
for some positive constant $K_1$ independent of $n$, and this is exactly \eqref{eq:uniformupperboundalpha}. 

 Letting $n\to+\infty$ in \eqref{eq:uniformupperboundalpha}, we then obtain
$$
u_{\alpha,p_0}(r)\leq  \frac{C}{\left(r^2- t_0^2 + K\right)^{\frac{\Nm-2}{2}}} \qquad \hbox{for} \ r \geq t_0\, .  
$$
 Since $p_0>p^*_->\frac{\tilde{N}_-}{\tilde{N}_--2}$, it  then follows that
 $$
\lim_{r\to +\infty}r^{\frac{2}{p_0-1}}u_{\alpha,p_0}(r)=0\, ,
$$
meaning that $u_{\alpha,p_0}$ is a  fast decaying solution of problem \eqref{eq:probMpm} with $p=p_0$. Since $\alpha<\alpha^*_-(p_0)$, this is again a contradiction to Theorem \ref{fastMm}. Hence, by the arbitrary choice of $\alpha<  \alpha^*_-(p_0)$, we deduce that $\liminf_{p\to p_0} \alpha^*_-(p)\geq \alpha^*_-(p_0)$, which finally proves the continuity of $\alpha^*_-$.
\end{proof}

\section{Critical exponents for the existence of radial sign-changing solutions in the ball}

In this section we prove Theorem \ref{mainteo1critexp}. Since the proof of ii) requires several steps we start by considering the case $\F=\Mp$, i.e. we consider the problem

\begin{equation}\label{eq:probM+p}
\begin{cases}
-\Mp(D^2 u)=|u|^{p-1}u & \hbox{in} \ B\\
\qquad  \ \ \ \ \ u=0 & \hbox{on} \ \partial B\\
\qquad  \ \ \ \ \ u(0)>0 &
\end{cases}
\end{equation}
Let us define the following set
\beq\label{defAnodalMp}
 \A:=\{ p \in (1,+\infty)\,: \ \text{there exists $u_p$ radial sign-changing solution to \eqref{eq:probM+p}} \}.
\eeq

\begin{remark}\label{rem1sectMP}
The set $\A$ is nonempty in view of \cite[Theorem 1.3]{GLP}. Moreover, by a trivial scaling argument, it is easy to check that $\A$ coincides with the set of $p \in (1,+\infty)$ for which there exists a nodal solution $u_p$ to \eqref{eq:probM+p} which changes sign exactly once.
\end{remark}
As a first result we show a crucial upper bound for $\sup \A$.
\begin{proposition}\label{prop1:sectexpcrtisc}
It holds that $$\sup \A < p^*_+.$$
\end{proposition}
\begin{proof}
We first observe that $\sup \A \leq p^*_+$, because for $p\geq\psp$ there cannot exist positive radial solutions to \eqref{eq:probM+p}. 

Now, assume by contradiction that $\sup \A = p^*_+$. Then we can find a sequence $(u_{p_n})_{p_n \in \A}$ of nodal radial solutions to \eqref{eq:probM+p}, with $p_n \nearrow \psp$. In view of Remark \ref{rem1sectMP} we can assume without loss of generality that $u_{p_n}$ changes sign exactly once. Let us consider the rescaled function
\beq\label{scaling1}
\tilde u_{p_n}(x)=r_1^{\frac{2}{p_n-1}} u_{p_n}(r_1 x), \;\;\; x \in B_{\frac{1}{r_1}},
\eeq
where $r_1=r_1(n)$ is the node of $u_{p_n}$.

By construction, for $x\in B$, $(\tilde u_{p_n})_{p_n \in \A}$ is a sequence of almost critical positive solutions of \eqref{eq:probM+p}. Then, in view of Proposition \ref{Prop:asanpos}, we have $\tilde u_{p_n}^\prime(1) \to 0$ as $p_n\to\psp$. In addition, by construction, the function $\tilde u_{p_n}^-(x):=r_1^{\frac{2}{p_n-1}} u^-_{p_n}(r_1 x)$, $x \in A_{1,\frac{1}{r_1}}$ is a positive radial solution of
\begin{equation*}
\begin{cases}
-\Mm(D^2 u)=u^{p_n} & \hbox{in} \ A_{1,\frac{1}{r_1}}\\
\qquad  \ \ \ \ \ u=0 & \hbox{on} \ \partial A_{1,\frac{1}{r_1}}.
\end{cases}
\end{equation*}
In particular, $\tilde u_{p_n}^-=\tilde u_{p_n}^-(r)$ satisfies \eqref{eq:initvalprob-} with $\alpha=\alpha(p_n)=(\tilde u^-_{p_n})^\prime(1) \to 0$, as $p_n\to\psp$, and we have $\tilde u_{p_n}^-(1/r_1)=0$.

On the other hand, since $\psp>\psm$ then in view of \cite[Theorem 1.1]{GIL} and Theorem \ref{fastMm} we have $\alpha_-^*(\psp)>0$ and, from Proposition \ref{lem:asalphastar}, we can find $\alpha_0>0$ such that $\alpha_-^*(p)>\alpha_0>0$ for all $p$ in a sufficiently small neighborhood of $\psp$. In particular, for any $\alpha \in (0,\alpha_0)$ and for any $p$ sufficiently close to $\psp$ the unique solution to \eqref{eq:initvalprob-} is defined and positive in the whole $(1,+\infty)$, but this contradicts the properties of $\tilde u^-_{p_n}$, namely that $\tilde u_{p_n}^-(1/r_1)=0$. 
This gives a contradiction and the proof is complete.
\end{proof}

\begin{proposition}\label{prop2:sectexpcrtisc}
It holds that $$\sup \A > \psm.$$
\end{proposition}
\begin{proof}
In order to prove the result we construct a sign-changing solution $u_p$ to \eqref{eq:probM+p} for $p$ in a sufficiently small right neighborhood of $\psm$.

To this end, let $p \in (\psm,\psp)$ and as in Sect. 2 we denote by $v_{p,+}$ the unique positive radial solution of
\begin{equation}\label{eq:probM+ppos}
\begin{cases}
-\Mp(D^2 u)=u^{p} & \hbox{in} \ B,\\
\qquad  \ \ \ \ \ u=0 & \hbox{on} \ \partial B.
\end{cases}
\end{equation}
Since $\psm<\psp$, choosing $0<\delta<(\psp-\psm)$ we have $I_\delta:=(\psm,\psm+\delta) \subset (\psm,\psp)$ and $m_p:=\max v_{p,+}=v_{p,+}(0)$ is uniformly bounded for $p \in I_\delta$. Moreover we can find $\gamma>0$ such that for all $p \in I_\delta$ it holds $|v_{p,+}^\prime(1)|>\gamma>0$. Hence, exploiting Proposition \ref{lem:asalphastar} and since $\alpha_-^*(p)\to 0$ as $p\to\psm$ (see the proof of Proposition \ref{lem:asalphastar}) we find a right neighborhood of $\psm$, let us say $I_\varepsilon:=(\psm,\psm+\varepsilon)$ with $\varepsilon<\delta$, such that  $|v_{p,+}^\prime(1)|>\alpha_-^*(p)$, for all $p \in I_\varepsilon$. This means that for all $p \in I_\varepsilon$,  if we take $\alpha=\alpha(p)= - v_{p,+}^\prime(1)$ in \eqref{eq:initvalprob-} then we have $\rho_\alpha=\rho_\alpha(p)<+\infty$ and the unique maximal positive solution $u_{\alpha(p)}=u_{\alpha(p)}(r)$ vanishes at $r=\rho_{\alpha(p)}$. Hence for $p \in I_\varepsilon$ we can glue the two solutions by defining $z_p:[0,\rho_{\alpha(p)}] \to \R$ as
$$z_p(r):=\begin{cases} v_{p,+}(r) & \text{if} \ r \in [0,1],\\ 
-u_{\alpha(p)}(r) & \text{if} \ r \in (1,\rho_{\alpha(p)}]. 
\end{cases}
$$
Then, setting $u_p(r):=[\rho_{\alpha(p)}]^{\frac{2}{p-1}}z_p(r\rho_{\alpha(p)})$ we easily check that $u_p(x)=u_p(|x|)$ is a radial sign-changing solution to \eqref{eq:probM+p}.
\end{proof}

By the very definition of $\sup \A$ we have that for $p>\sup \A$ radial sign-changing solution to \eqref{eq:probM+p} cannot exist.
 In the next proposition we show that the same happens for $p=\sup \A$. 
 
 \begin{proposition}\label{prop:nonexistscsolMP}
It holds that  $\sup \A\notin\A$.
\end{proposition}
\begin{proof}
Assume by contradiction that there exists a a radial sign-changing solution $u_{*}$ to \eqref{eq:probM+p} for $p=\sup \A$.  In view of Remark \ref{rem1sectMP} we can assume without loss of generality that $u_{*}$ changes sign exactly once and we denote by $r_{*} \in (0,1)$ its nodal radius. 

Let us consider the rescaled radial function $\tilde u_{*}(x):=r_{*}^{\frac{2}{p-1}} u_{*}(r_{*}x)$, $x \in B_{1/r_{*}}$ and set $\alpha_*:= -\tilde u_{*}^\prime(1)>0$. Then, by construction $(\tilde u_*)^-(r)$ (the negative part of $\tilde u_*$, defined by taking the maximum between $-\tilde u_*$ and zero) coincides for $r \in [1,1/r_*]$ with the unique maximal positive solution to \eqref{eq:initvalprob-} with $\alpha=\alpha_*$, $p=\sup \A$. 

Therefore, since $\sup \A>\psm$ and $(\tilde u_*)^-(1/r_{*})=0$, from Theorem \ref{fastMm} we have
\beq\label{eq:stimaalphastar}
\alpha_*>\alpha_-^*(\sup \A).
\eeq
Now we observe that, since $\sup \A<\psp$, up to a subsequence, as $p_n\searrow\sup \A$ the unique positive  solution $v_{p_n,+}$ to \eqref{eq:probM+ppos} converges in $C^2(\overline{B})$ to $\tilde u_{*}\big|_{\overline B}$. In particular this implies that $-v_{p_n,+}^\prime(1) \to \alpha_*$, for some sequence $p_n\searrow \sup \A$. Moreover, from Proposition \ref{lem:asalphastar} and \eqref{eq:stimaalphastar} we infer that $-v_{p_n,+}^\prime(1)> \alpha_-^*(p_n)$ for all $p_n$ sufficiently close to $\sup \A$.  

Therefore, fixing $p_n>\sup \A$ sufficiently close to $\sup \A$, and taking $\alpha=\alpha(p_n)= -v_{p_n,+}^\prime(1)>0$ the unique maximal positive solution $u_{\alpha(p_n)}$ to \eqref{eq:initvalprob-} vanishes at some radius $\rho_{\alpha(p_n)}$ such that $1<\rho_{\alpha(p_n)}<+\infty$. Then, gluing $v_{p_n,+}$ and $u_{\alpha(p_n)}$ as in the proof of Proposition \ref{prop2:sectexpcrtisc} we obtain a radial sign-changing solution to Problem \eqref{eq:probM+p}. Hence $p_n \in \A$ which is a contradiction since $p_n>\sup \A$. The proof is complete.
\end{proof}

\begin{proposition}\label{lemmacriticalexpkzones}
For any integer $k\geq 2$ and any $p \in \A$ there exists a radial sign-changing solution $u_p$ to \eqref{eq:probM+p} with $k$ nodal regions.
\end{proposition}
\begin{proof}
Let $p\in \A$. We argue by induction on $k$. 
The basic step $k=2$ is obvious by the definition of $\A$ (see also Remark \ref{rem1sectMP}).

Assume that there exists $u_{p,k}$ radial sign-changing solution to \eqref{eq:probM+p} with $k$ nodal domains. We need to distinguish between two cases. 

If $k$ is even, then $u_{p,k}<0$ in the last nodal region and thus by Hopf's Lemma we infer that $u_{p,k}^\prime(1)>0$. 
Then, since $p<\sup \A<\psp$, from \cite[Theorem 1.1]{GIL} we have that for any $\alpha>0$ the unique maximal positive solution $u_{\alpha}=u(\alpha,p,r)$ of the initial value problem \eqref{eq:initvalprob+} vanishes at some $\rho_\alpha<+\infty$. Hence, taking $\alpha=\alpha(p)=u_{p,k}^\prime(1)>0$,  gluing $u_{p,k}$, $u_{\alpha(p)}$ as in the proof of Proposition \ref{prop2:sectexpcrtisc} and rescaling, we obtain a radial solution $u_{p,k+1}$ of \eqref{eq:probM+p} with $k+1$ nodal regions and such that $u_{p,k+1}(0)>0$. This complete the proof of the inductive step when $k$ is even.

If $k$ is odd the previous argument works only for $p\leq\psm$ (see  \cite[Theorem 1.2]{GIL}). For this reason we proceed in a different way. Let $u_{p,k}$ be a radial sign-changing solution to \eqref{eq:probM+p} with $k$ nodal regions. Since $k$ is odd then $u_{p,k}>0$ in the last nodal region and thus $-u_{p,k}^\prime(1)>0$. We claim that $-u_{p,k}^\prime(1)>\alpha_-^*(p)$, where $\alpha_-^*$ is the critical slope defined in \eqref{defalphastar}. We first observe that if the claim is true then the maximal positive solution $u_{\alpha}=u(\alpha,p,r)$ of \eqref{eq:initvalprob-} with $\alpha=\alpha(p)=-u_{p,k}^\prime(1)$ vanishes at some $\rho_{\alpha(p)}<+\infty$, and thus we can construct $u_{p,k+1}$ satisfying the desired properties by gluing $u_{p,k}$ and $u_{\alpha(p)}$ in the same way as before. 

To prove the claim, let $v_{p,+}$ be the only  positive solution to \eqref{eq:probM+ppos} and consider its trajectory $\gamma_1(t)=(x_1(t), x^\prime_1(t))$, $t \in (-\infty, 0]$ in the phase-plane, where $x_1$ is the Emden-Fowler transform of $v_{p,+}$, obtained by setting $r=e^t$ and
$$
x_1(t):=e^{\frac{2}{p-1}t} v_{p,+}(e^t).
$$
By construction $\gamma_1$ lies in the right-half plane and it is elementary to check that $\gamma_1(t)\to (0,0)$ as $t\to -\infty$, and $\gamma_1(0)=(0,v_{p,+}^\prime(1))$, with $v_{p,+}^\prime(1)<0$. On the other hand, if we transform the restriction of $u_{p,k}$ to its last nodal component we obtain a trajectory $\gamma_2(t)=(x_2(t),x_2^\prime(t))$, $t \in [\log(r_{k-1}),0]$ lying in the right-half plane and such that $\gamma_2(\log(r_{k-1}))=(0,r_{k-1}^{\frac{p+1}{p-1}}u_{p,k}^\prime(r_{k-1}))$ with $u_{p,k}^\prime(r_{k-1})>0$, while $\gamma_2(0)=(0,u_{p,k}^\prime(1))$ and $u_{p,k}^\prime(1)<0$.  Moreover, since $x_1$, $x_2$ satisfy the same autonomous ODE (see \cite{GIL} for more details) then $\gamma_1$ and $\gamma_2$ cannot intersect and thus it follows that $0>v_{p,+}^\prime(1)>u_{p,k}^\prime(1)$. 

Now, let us prove that $-v'_{p,+}(1)>\alpha_-^*(p)$. Indeed, since $1<p<\sup \A$ and $p\in \A$ we know that there exists a radial sign-changing $u_{p,2}$ solution to \eqref{eq:probM+p} with two nodal regions. Then, denoting by $r_1\in(0,1)$ the node of $u_{p,2}$ and considering the usual scaling $\tilde u_{p,2}(x):=r_1^{\frac{2}{p-1}} u_{p,2}(r_1 x)$, $x \in B_{1/r_1}$, we infer that the restriction $\tilde u_{p,2}\big|_B$ coincides with $v_{p,+}$ (uniqueness of the positive radial solution) and the restriction $\tilde u_{p,2}^-\big|_{A_{1,1/r_1}}$ coincides with the unique positive solution to
\eqref{eq:initvalprob-} with $\alpha=-v_{p,+}^\prime(1)$. Therefore, since $\tilde u_{p,2}\big|_{A_{1,1/r_1}}$ vanishes at $1/r_1$ then by definition of $\alpha_-^*$ and by Theorem \ref{fastMm} we conclude that $-v_{p,+}^\prime(1)>\alpha_-^*(p)$.

Then we have proved that $\alpha_-^*(p)<-v_{p,+}^\prime(1)<-u_{p,k}^\prime(1)$. This concludes the proof.
\end{proof}
\begin{proof}[Proof of Theorem \ref{mainteo1critexp}]
To prove i) we observe that the existence of radial sign-changing solutions to \eqref{eq:probMgen} for $\F=\Mm$  is a consequence of the existence of the positive solution of the same problem, combined with the  Liouville type results obtained in \cite{GIL}. For this, let $p<\psm$ and let $v_p$ be the  positive solution of 
\begin{equation}\label{eqnuovascmm}
\begin{cases}
-\Mm(D^2 u)=|u|^{p-1}u & \hbox{in} \ B\\
\qquad  \ \ \ \ \ u=0 & \hbox{on} \ \partial B\\
\qquad  \ \ \ \  u(0)>0 &
\end{cases}
\end{equation}
By Hopf's Lemma it holds that $|v_p^\prime(1)|>0$. Let $w_p$ be the positive maximal solution of \eqref{eq:initvalprob+} with initial slope $w_p'(1)=|v_p^\prime(1)|$. Since $p<\psm$ (which in particular implies that $p<\psp$) then in view of \cite[Theorem 3.1]{GIL}, the function $w_p$ must vanish at some $\rho=\rho(p)>1$. Hence the function
\begin{equation*}
u_p(r)=\begin{cases}
v_p(r) & \text{if $r\leq1$}\\
-w_p(r) & \text{if $r\in(1,\rho]$}
\end{cases}
\end{equation*}
defines a radial sign-changing solution, with two nodal regions, of
$$
-\Mm(D^2u)=|u|^{p-1}u\qquad\text{in $B_{\rho}$}
$$
and such that $u_p(0)>0$. By scaling,   $\tilde u_p(r)=\rho^{\frac{2}{p-1}}u_p(\rho r)$ is a sign-changing  solution of \eqref{eqnuovascmm}. This completes the proof of i) in the case of two nodal regions. 

Let us point out that such gluing argument can be performed inductively, so providing the existence of sign-changing solutions for any number of nodal regions. The key point in this procedure is that $p$ is subcritical both for ${\mathcal M}^-_{\lambda,\Lambda}$ and ${\mathcal M}^+_{\lambda,\Lambda}$, which implies that for any choice of the initial slope $\alpha>0$ the unique positive solution of \eqref{eq:initvalprob-} or \eqref{eq:initvalprob+} vanishes at some finite $\rho \in (1,+\infty)$ (see \cite[Theorem 3.1]{GIL}).\\

Let us prove ii). We define
 \beq\label{defespcritscMp}
p^{**}_+:=\sup \A,
\eeq
where $\A$ is given by \eqref{defAnodalMp}. Then ii) is a consequence of Proposition \ref{prop1:sectexpcrtisc}, Proposition \ref{prop2:sectexpcrtisc}, Proposition \ref{prop:nonexistscsolMP} and Proposition \ref{lemmacriticalexpkzones}.
\end{proof}

\section{Asymptotic analysis of radial sign-changing solutions to \eqref{eq:probM-} with two nodal regions}\label{Sec3}
\noindent
In this section $u_\e$  will denote a radial sign-changing solution of \eqref{eq:probM-} with two nodal regions. We set 
\begin{equation}
M_0=M_0(\e):=\left\|u^+_\e\right\|_\infty=u^+_\e(0),
\end{equation}
where $u_\e^+:=\max\{0,u_\e\}$ is the positive part of $u_\e$, and we denote by $r_1=r_1(\e) \in (0,1)$ the node of $u_\e$, i.e. the unique point $r_1 \in (0,1)$ such that $u_\e(r_1)=0$. As in the previous section, $v_{p_\e,-}$ stands for the only positive  solution of \eqref{eq:probM-}. We begin with a preliminary result.

\begin{proposition}\label{scaling}
The following statements hold:
\begin{itemize}
	\item[i)] $M_0>\left\|v_{p_\e,-}\right\|_\infty$;
	\item[ii)] $\displaystyle r_1=\left(\frac{\left\|v_{p_\e,-}\right\|_\infty}{M_0}\right)^{\frac{p_\e-1}{2}}$;
	\item[iii)] $\displaystyle \lim_{\e\to0}M_0^{\frac{p_\e-1}{2}}r_1=+\infty$;
	\item[iv)] $\displaystyle r_1^{\frac{p_\e+1}{p_\e-1}}u_\e^\prime(r_1)=(v_{p_\e,-})^\prime(1)$.
\end{itemize}
\end{proposition}
\begin{proof}
 Consider the rescaled function $\tilde u_\e(x):=r_1^{\frac{2}{p_\e-1}} u_\e(r_1 x)$, $x\in B$. It is elementary to check that $\tilde u_\e$ is a positive radial solution to \eqref{eq:probM-}. Hence, by uniqueness, we infer that $\tilde u_\e=v_{p_\e,-}$, and the result easily follows from Proposition \ref{Prop:asanpos}.
\end{proof}

Concerning the negative part of $u_\e$, namely  $u^-_\e:=\max\left\{-u_\e,0\right\}$, we adopt the following notations: $$M_1=M_1(\e):=\left\|u^-_\e\right\|_\infty,$$ $s_1=s_1(\e)\in(r_1,1)$ is the point where the maximum $M_1$ is attained, i.e. $M_1=u^-_\e(s_1)$, and $t_1\in(s_1,1)$ is the only radius such that 
\begin{equation*}
(u_\e^-)''(r)<0\quad\text{for $r\in(r_1,t_1)$}, \ \ \ \ \ (u_\e^-)''(r)>0\quad\text{for $r\in(t_1,1)$.}
\end{equation*}
\begin{remark}\label{rem:firsteigen}
We point out that $M_1$ is bounded away from zero, in fact it is bounded from below by  the principal eigenvalue $\lambda^+_1=\lambda^+_1(-\Mp;B)$. Indeed $u^-_\e$  satisfies in the annulus $A_{r_1,1}=\{x\in \RN; \ r_1 <|x|<1\}$ the following
\begin{equation*}
\left\{\begin{array}{cl}
-\Mp(D^2u^-_\e)\leq M_1^{p_\e-1}u^-_\e & \text{in $A_{r_1,1}$}\\
u^-_\e=0 & \text{on $\partial A_{r_1,1}$}.
\end{array}\right.
\end{equation*}
Since the principal eigenvalue $\lambda^+_1$ gives a threshold for the validity of the maximum principle (see \cite{BEQ}) and $u^-_\e>0$ in $A_{r_1,1}$, then necessarily
\begin{equation}\label{bound M1}
M_1^{p_\e-1}\geq \lambda^+_1(-\Mp;A_{r_1,1})\geq \lambda^+_1(-\Mp;B).
\end{equation}
\end{remark}

\begin{proposition}\label{cor1}
The following statements hold:
\begin{itemize}
	\item[i)] $\displaystyle\lim_{\e\to0}r_1=0$;
	\item[ii)] $\displaystyle \lim_{\e\to0} M_1^{\frac{p_\e-1}{2}}r_1=0$;
	\item[iii)] $\displaystyle \lim_{\e\to0} \frac{r_1}{s_1}=0$;
	\item[iv)]  $\displaystyle \lim_{\e\to0} |u_\e'(r_1)|=+\infty$.
\end{itemize}
\end{proposition}
\begin{proof}
Let us consider the energy-like functional 
\begin{equation}\label{E1}
H(r):=\frac{(u^\prime_\e(r))^2}{2}+\frac{|u_\e(r)|^{p_\e+1}}{\lambda(p_\e+1)}, \ \ r\in[0,1].
\end{equation}
By a straightforward computation it holds that $H^\prime\leq0$ in $[r_1,s_1]$ (actually it could be proved that  $H^\prime\leq0$ in  the whole interval $[0,1]$). This in particular yields
\begin{equation}\label{bound E1}
\frac{M_1^{p_\e+1}}{\lambda(p_\e+1)}=H(s_1)\leq H(r_1)=\frac{(u^\prime_\e(r_1))^2}{2}\,.
\end{equation}
Using  Proposition \ref{scaling}-iv), Proposition \ref{Prop:asanpos}-i,vi),  \eqref{bound M1} and \eqref{bound E1} we then obtain
\begin{equation}\label{8/10eq3}
\begin{split}
r_1^{2\frac{p_\e+1}{p_\e-1}}&=\frac{(v^\prime_{p_\e,-}(1))^2}{(u^\prime_\e(r_1))^2}
\leq\frac{\lambda(p_\e+1)}{2}\frac{(v^\prime_{p_\e,-}(1))^2}{M_1^{p_\e+1}}\\
&\leq\frac{\lambda(p_\e+1)}{2}\frac{(v^\prime_{p_\e,-}(1))^2}{\lambda^+_1(-\Mp,B)^\frac{p_\e+1}{p_\e-1}}\to0\quad\text{as $\e\to0$}.
\end{split}
\end{equation}
This proves i).  Moreover, again by \eqref{bound E1} and Proposition \ref{scaling}-iv)
$$
r_1 M_1^{\frac{p_\e-1}{2}}\leq\left(\frac{\lambda(p_\e+1)}{2}\right)^{\frac{p_\e-1}{2(p_\e+1)}}r_1 |u^\prime_\e(r_1)|^{\frac{p_\e-1}{p_\e+1}}=\left(\frac{\lambda(p_\e+1)}{2}\right)^{\frac{p_\e-1}{2(p_\e+1)}}|v^\prime_{p_\e,-}(1)|^{\frac{p_\e-1}{p_\e+1}},
$$
which proves ii). 

\medskip
In order to prove iii) and iv) let us consider the energy-like functionals
\begin{equation}\label{E2}
H_{\gamma,\eta}(r):=r^{\Nm}\left(\frac{(u^\prime_\e(r))^2}{2}+\frac{\gamma}{p_\e+1}|u_\e(r)|^{p_\e+1}\right)+\eta r^{\Nm-1}u_\e(r)u^\prime_\e(r),
\end{equation}
where $\gamma, \eta$ are real parameters to be chosen later.
It is easy to check that  for $r\in[r_1,s_1]$ 
\begin{equation}\label{derE2}
\begin{split}
H^\prime_{\gamma,\eta}(r)&=\left(\eta+1-\frac\Nm2\right)r^{\Nm-1}(u^\prime_\e(r))^2\\&\quad+\left(\frac{\gamma\Nm}{p_\e+1}-\frac\eta\lambda\right)r^{\Nm-1}|u_\e(r)|^{p_\e+1}\\&\qquad+\left(\gamma-\frac1\lambda\right)r^{\Nm}|u_\e(r)|^{p_\e-1}u_\e(r)u^\prime_\e(r).
\end{split}
\end{equation}
Choose $$\eta=\frac{\Nm-2}{2},\quad\gamma=\frac{\eta(p_\e+1)}{\lambda\Nm}$$
in \eqref{derE2}. Since $p_\e>\frac{\Nm+2}{\Nm-2}$ for small $\e>0$, then  $H^\prime_{\gamma,\eta}\geq0$ in $[r_1,s_1]$. Hence $H_{\gamma,\eta}(r_1)\leq H_{\gamma,\eta}(s_1)$ which reads as 
\begin{equation}\label{eq1}
r_1^{\Nm}(u^\prime_\e(r_1))^2\leq\frac{\Nm-2}{\lambda\Nm}s_1^{\Nm}M_1^{p_\e+1}.
\end{equation} 
By the convexity of $u_\e$ in $[r_1,s_1]$ we also have
\begin{equation}\label{eq2}
M_1^2\leq (u^\prime_\e(r_1))^2s_1^2.
\end{equation}
Putting together \eqref{eq1}-\eqref{eq2} we get
$$
\left(\frac{r_1}{s_1}\right)^{\Nm}\leq\frac{\Nm-2}{\lambda\Nm}M_1^{p_\e-1}s_1^2=\frac{\Nm-2}{\lambda\Nm}M_1^{p_\e-1}r_1^2\left(\frac{s_1}{r_1}\right)^2.
$$
Hence, by ii), iii) follows. 

\medskip

Now set in \eqref{E2}-\eqref{derE2}
$$\eta=\frac{\Nm-2}{2},\quad\gamma=\frac{1}{\lambda}.$$
 With such a choice we easily check that $H_{\gamma,\eta}^\prime\leq0$ in $[r_1,s_1]$. Then $H_{\gamma,\eta}(r_1)\geq H_{\gamma,\eta}(s_1)$ and  using \eqref{bound M1} we infer that
\begin{equation}\label{equtile}
\begin{split}
(u^\prime_\e(r_1))^2&\geq\frac{2}{\lambda(p_\e+1)}\left(\frac{s_1}{r_1}\right)^{\Nm} M_1^{p_\e+1}\\
&\geq\frac{2}{\lambda(p_\e+1)}\left(\frac{s_1}{r_1}\right)^{\Nm} \left(\lambda^+_1(-\Mp,B)\right)^{\frac{p_\e+1}{p_\e-1}}.
\end{split}
\end{equation}
This implies iv) because of iii). 
\end{proof}
In the following statements we will make use of the rescaled function defined by
\beq\label{defrescaledtwonodal}
\hat u^-_\e(x):=\frac{1}{M_1}u^-_\e\left(\frac{x}{M_1^{\frac{p_\e-1}{2}}}\right)\qquad x\in\hat A_\e,
\eeq
where $\hat A_\e$ is the annulus $\hat A_\e:=A_{r_1 M_1^{\frac{p_\e-1}{2}},M_1^{\frac{p_\e-1}{2}}}$. 
The function $\hat u^-_\e$ is a positive radial solution of
\begin{equation}\label{eq:defrescaledtwo}
\left\{\begin{array}{ll}
-\Mp(D^2 u)=u^{p_\e} & \text{in $\hat A_\e$},\\[4pt]
 \ \ \ \  \ \ \ \  \ \ \ \ \  \ \  \ u=0& \text{on $\partial \hat A_\e$}. 
\end{array}\right.
\end{equation}
For convenience of notations we set
\begin{equation}\label{notation}
\hat r_1=\hat r_1(\e):=r_1 M_1^{\frac{p_\e-1}{2}},\quad \hat s_1=\hat s_1(\e):=s_1 M_1^{\frac{p_\e-1}{2}}, \quad \hat t_1=\hat t_1(\e):=t_1 M_1^{\frac{p_\e-1}{2}}.
\end{equation}
The first result is about the asymptotic behavior of $\hat s_1$.
\begin{proposition}\label{prop sigma}
 $\displaystyle\lim _{\e\to 0}\hat s_1= 0$.
\end{proposition}
\begin{proof}
We first prove that $\hat s_1$ is bounded from above. For this let us consider the energy  function \eqref{E1}, which is monotone decreasing in $[r_1,s_1]$ (actually in $[0,1]$). Hence for any $r\in[r_1,s_1]$ we have $H(r)\geq H(s_1)$ and 
$$
(u^-_\e)^\prime(r)\geq\sqrt{\frac{2}{\lambda(p_\e+1)}}\sqrt{\left(M_1^{p_\e+1}-(u^-_\e(r))^{p_\e+1}\right)}.
$$
Integrating in $[r_1,s_1]$ we infer that
\begin{equation}\label{eq3}
\int_{r_1}^{s_1}\frac{(u^-_\e)'(r)}{\sqrt{\left(M_1^{p_\e+1}-u^-_\e(r)^{p_\e+1}\right)}}\,dr\geq\sqrt{\frac{2}{\lambda(p_\e+1)}}(s_1-r_1).
\end{equation}
With the change of variable $t=\frac{u^-_\e}{M_1}$ we obtain
$$
\int_{r_1}^{s_1}\frac{(u^-_\e)'(r)}{\sqrt{\left(M_1^{p_\e+1}-u^-_\e(r)^{p_\e+1}\right)}}\,dr=\frac{1}{M_1^{\frac{p_\e-1}{2}}}\int_0^1\frac{1}{\sqrt{1-t^{p_\e+1}}}\,dt.
$$
Then from \eqref{eq3} we deduce that
$$
\hat s_1-\hat r_1\leq\sqrt{\frac{\lambda(p_\e+1)}{2}}\int_0^1\frac{dt}{\sqrt{1-t^{p_\e+1}}}.
$$
Sending $\e\to0$ and using  Proposition \ref{cor1}-ii) 
$$
\limsup_{\e\to0}\hat s_1\leq\sqrt{\frac{\lambda(p^*_-+1)}{2}}\int_0^1\frac{dt}{\sqrt{1-t^{p^*_-+1}}}\leq\frac\pi2\sqrt{\frac{\lambda N}{N-2}},
$$
which proves the claim. 

 From \eqref{case1} it is easy to check that
\begin{equation}\label{eq4*}
\left(r^{\Nm-1}( u^-_\e)^\prime\right)^\prime+\frac{1}{\lambda}r^{\Nm-1} u^-_\e(r)^{p_\e}=0\qquad\text{in\; $[r_1, s_1]$.}
\end{equation}
Integrating from $r_1$ to $r\in[r_1,s_1]$ we obtain
$$
-r^{\Nm-1}(u^-_\e)^\prime(r)+r_1^{\Nm-1}(u^-_\e)^\prime(r_1)=\frac1\lambda\int_{r_1}^rs^{\Nm-1}u^-_\e(s)^{p_\e}\,ds\leq\frac{M_1^{p_\e}}{\lambda\Nm}r^{\Nm}
$$
and then
$$
-(u^-_\e)^\prime(r)\leq\frac{M_1^{p_\e}}{\lambda\Nm}r-r_1^{\Nm-1}(u^-_\e)^\prime(r_1)r^{1-\Nm}\qquad\text{in\; $[r_1,s_1]$.}
$$
Integrating the above inequality in $[r_1,s_1]$, we infer that 
\begin{equation}\label{eq1-2/10}
\frac{r_1( u^-_\e)^\prime(r_1)}{\Nm-2}\left[1-\left(\frac{r_1}{s_1}\right)^{\Nm-2}\right]\leq M_1\left(1+\frac{\hat s_1^2}{2\lambda\Nm}\right).
\end{equation}
Since $\hat s_1$ is bounded from above and $\displaystyle\lim_{\varepsilon\to0}\frac{r_1}{s_1}\to0$, in view of Proposition \ref{cor1}-iii),  then 
$$ r_1(u^-_\e)^\prime(\hat r_1)\leq CM_1$$ for some positive constant $C$. Moreover $r_1(u^-_\e)^\prime(r_1)\geq0$. Hence using  \eqref{equtile}
\begin{equation*}
\begin{split}
C^2&\geq \left(r_1\frac{(u^-_\e)^\prime(r_1)}{M_1}\right)^2\\
&\geq\frac{2}{\lambda(p_\e+1)}\left(\frac{s_1}{r_1}\right)^{\Nm-2} \hat s^2_1.
\end{split}
\end{equation*}
The conclusion follows by Proposition \ref{cor1}-iii).
\end{proof}

From the previous result we immediately deduce that $\hat r_1\to 0$, as $\e\to 0$. Moreover we have 
\begin{corollary}\label{cors1zero}
$\displaystyle \lim_{\e\to0}s_1=0$.
\end{corollary}
\begin{proof}
Use \eqref{bound M1} and Proposition \ref{prop sigma}.
\end{proof}
Next we study the asymptotic behavior of $\hat t_1$ and show that it cannot converge to zero. We begin with a stronger result.
\begin{proposition}\label{tau}
 $\displaystyle \liminf_{\e\to0}t_1 (u^-_\e(t_1))^{\frac{p_\e-1}{2}}\geq\sqrt{\frac{2\lambda(N-1)}{p^*_-+1}}$.
\end{proposition}
\begin{proof}
Let us consider the energy-like functional
$$H(r):=r^{N}\left(\frac{[(u^-_\e)^\prime(r)]^2}{2}+\frac{1}{\lambda(p_\e+1)}(u^-_\e(r))^{p_\e+1}\right)+\frac{N}{p_\e+1} r^{N-1}u^-_\e(r)(u^-_\e)^\prime(r), \ \ r \in [s_1,t_1].$$
By \eqref{case2}, a direct computation shows that $H^\prime(r)\geq0$. Hence $H(s_1) \leq H(t_1)$ and thus
$$
s_1^N\frac{M_1^{p_\e+1}}{\lambda(p_\e+1)}\leq t_1^N\left(\frac{[(u^-_\e)^\prime(t_1)]^2}{2}+\frac{1}{\lambda(p_\e+1)}(u^-_\e(t_1))^{p_\e+1}\right)+\frac{N}{p_\e+1} t_1^{N-1}u^-_\e(t_1)(u^-_\e)^\prime(t_1).
$$
Exploiting the equation
$$
(N-1)\frac{(u^-_\e)'(t_1)}{t_1}=-\frac1\lambda (u^-_\e(t_1))^{p_\e},
$$ 
we have
 $$
s_1^N\frac{M_1^{p_\e+1}}{\lambda(p_\e+1)}\leq\frac{1}{\lambda(N-1)}t_1^N(u^-_\e(t_1))^{p_\e+1}\left(\frac{1}{2\lambda(N-1)}t_1^2(u^-_\e(t_1))^{p_\e-1}-\frac{1}{p_\e+1}\right).
$$
This implies that
$$
t_1^2(u^-_\e(t_1))^{p_\e-1}\geq\frac{2\lambda(N-1)}{p_\e+1}
$$
and the conclusion follows passing to the limit as $\e\to0$.
\end{proof}
As an immediate consequence of the previous proposition we get
\begin{corollary}\label{cort}
$\displaystyle\liminf_{\e\to 0} \hat t_1 \geq\sqrt{\frac{2\lambda(N-1)}{p^*_-+1}}.$
\end{corollary}
\begin{proof}
It suffices to observe that $\hat t_1 = t_1 M_1^{\frac{p_\e-1}{2}} \geq t_1 (u^-_\e(t_1))^{\frac{p_\e-1}{2}}$ so that the assertion follows by Proposition \ref{tau}.
\end{proof}
Coming back to the study of $u^-_\e$, from the previous results we obtain
\begin{proposition}\label{prop6Sect4}
$\displaystyle\limsup_{\e\to0}M_1<+\infty$.
\end{proposition}
\begin{proof}
Arguing by contradiction let us assume that along some subsequence, still denoted by $\e$, $M_1\to\infty$. Consider the rescaled function $\hat u^-_\e$ defined in \eqref{defrescaledtwonodal}. By construction  $\hat u^-_\e$ satisfies

$$
\hat u^-_\e(\hat s_1)=1,\quad 0\leq\hat u^-_\e\leq1.
$$

The limit of the domains $\hat A_\e$ is $\RN\backslash\left\{0\right\}$ and, since $(\hat u^-_\e)_\e$ is uniformly bounded and solves \eqref{eq:defrescaledtwo}, by regularity  estimates, we have that $\hat u^-_\e\to \hat u$ in $C^2_{{\rm loc}}(\RN\backslash\left\{0\right\})$, for some radially symmetric function $\hat u$, $0\leq \hat u\leq1$, which solves
\begin{equation}\label{eq7}
-\Mp(D^2 u)= u^{p^*_-}\qquad\text{in $\RN\backslash\left\{0\right\}$}.
\end{equation}
Now we want to show that $\hat u$ can be extended to a $C^2$ solution of \eqref{eq7} in the whole $\RN$ with $\hat u(0)=1$ and $\hat u^\prime(0)=0$. In the interval $(\hat s_1,\hat t_1)$ we have that $(\hat u^-_\e)^\prime\leq0$ and $(\hat u^-_\e)^{\prime\prime}\leq0$.
Therefore the equation satisfied by $\hat u^-_\e$ is 
$$
-(\hat u^-_\e)^{\prime\prime}(r)-\frac{N-1}{r}(\hat u^-_\e)^\prime(r)=\frac{(\hat u^-_\e)(r)^{p_\e}}{\lambda}\qquad r\in(\hat s_1,\hat t_1)
$$
that we can write as:
\begin{equation}\label{eq8}
\left(r^{N-1}(\hat u^-_\e)^\prime\right)^\prime=-\frac{1}{\lambda}r^{N-1}(\hat u^-_\e)(r)^{p_\e}\geq-\frac{1}{\lambda}r^{N-1}\,.
\end{equation}
Integrating between $\hat s_1$ and $r\in(\hat s_1,\hat t_1)$ we get
$$
r^{N-1}(\hat u^-_\e)^\prime(r)-\hat s_1^{N-1}(\hat u^-_\e)^\prime(\hat s_1)\geq-\frac{1}{\lambda N}(r^N-\hat s_1^N),
$$
and, since $(\hat u^-_\e)^\prime(\hat s_1)=0$ and $\hat s_1>0$, we obtain
\begin{equation}\label{delta}
(\hat u^-_\e)^\prime(r)\geq-\frac{1}{\lambda N}r\qquad \text{for $r\in(\hat s_1,\hat t_1)$}.
\end{equation}
Integrating again between $\hat s_1$ and $r\in(\hat s_1,\hat t_1)$, taking into account that $\hat u^-_\e(\hat s_1)=1$, we have
\begin{equation}\label{eq9}
\hat u^-_\e(r)\geq1-\frac{1}{2\lambda N}r^2\qquad \text{for $r\in(\hat s_1,\hat t_1)$}.
\end{equation}
Passing to the limit as $\e\to0$, taking into account  Proposition \ref{prop sigma} and  Corollary \ref{cort}, we infer that
$$
\hat u(r)\geq1-\frac{1}{2\lambda N}r^2\qquad \text{for $r\in\left(0,\sqrt{\frac{2\lambda(N-1)}{p^*_-+1}}\right)$}.
$$
From this, since $\hat u\leq1$, we deduce that
\begin{equation}\label{eq10}
\lim_{r\to0}\hat u(r)=1.
\end{equation}
Hence $\hat u$ can be extended by continuity at the origin, by setting $\hat u(0):=1$. Next we show that also $\hat u^\prime$ can be extended by continuity in $0$. From \eqref{delta} we have
$$
|(\hat u^-_\e)^\prime(r)|\leq\frac{r}{\lambda N}\qquad\text{for any $r\in(\hat s_1,\hat t_1)$}
$$
and passing to the limit as $\e\to0$, we get
\begin{equation}\label{eq11}
|\hat u^\prime(r)|\leq\frac{r}{\lambda N}\qquad\text{for $r\in\left(0,\sqrt{\frac{2\lambda(N-1)}{p^*_-+1}}\right)$}
\end{equation}
which gives $\lim_{r\to0}\hat u^\prime(r)=0$. Hence $\hat u$ is a positive radial solution of \eqref{eq7} that extends to a $C^1$ function near the origin. This implies that $\hat u$ is a (positive) $C^2$ radial solution of
\beq\label{eqhop2}
-\Mp(D^2 u)= u^{p^*_-} \ \ \hbox{in} \ \RN.
\eeq
Indeed, since $\hat u$ verifies $(r^{N-1} \hat u^\prime(r))^\prime=-\frac{1}{\lambda} r^{N-1} (\hat u(r))^{\psm}$ for $r\in\left(0,\sqrt{\frac{2\lambda(N-1)}{p^*_-+1}}\right)$, then, fixing $0<\delta<r$ and integrating between $\delta$ and $r$, we get that
$$r^{N-1} \hat u^\prime(r) - \delta^{N-1} \hat u^\prime(\delta)=-\frac{1}{\lambda} \int_\delta^r s^{N-1} (\hat u(s))^{\psm} \ ds.$$
Passing to the limit as $\delta \to 0$, by \eqref{eq11}, we obtain
\beq\label{eqhopital}
 \frac{\hat u^\prime(r)}{r}= -\frac{1}{\lambda r^N} \int_0^r s^{N-1} (\hat u(s))^{\psm} \ ds.
\eeq
By de L'$\rm{H\hat opital}$'s rule, the right-hand side of \eqref{eqhopital} has a finite limit as $r\to 0$ and thus the same holds for the left-hand side, which readily implies that $\hat u$ extends to a $C^2$ radial solution of \eqref{eqhop2}.

At the end, since $p^*_-<p^*_+$, then by \cite[Theorem 1.1]{FQ} we know that \eqref{eqhop2} has only the trivial solution, which contradicts the positivity of $\hat u$.
\end{proof}

Summing up, we have all the ingredients to prove the following
\begin{theorem}\label{teo1Sect4}
Let $u_\e$ be a radial sign-changing solution to \eqref{eq:probM-} with two nodal regions. Then, up to a subsequence, as $\e \to 0^+$ we have that $u_\e \to \bar u$ in $C^2_{loc}(\overline B\setminus\{0\})$, where $\bar u$ is the unique  negative solution of \eqref{eq:probM-crit2}
\end{theorem}
\begin{proof}
Let us consider the restriction of $u^-_\e$ to the annulus $A_{r_1,1}$. From Proposition \ref{prop6Sect4} we have that $u^-_\e\big|_{A_{r_1,1}}$ is uniformly bounded and from Proposition \ref{cor1}, i) we have $r_1\to 0$. Hence, by standard regularity theory, up to a subsequence as $\e\to 0$, we get that $u^-_\e \to \bar u^-$ in $C^2_{loc}(\overline{B}\setminus\{0\})$, where $\bar u^-$ is a non-negative radial solution of 
\beq\label{eqprobwholeball}
\begin{cases}
-\Mp(D^2 u)=|u|^{p_-^*-1}u & \hbox{in} \ B\setminus\{0\},\\
\qquad  \ \ \ \ \ u=0 & \hbox{on} \ \partial B.
\end{cases}
\eeq
We claim that $\bar u^-$ can be extended to a smooth positive solution of \eqref{eqprobwholeball} in the whole ball. For this, taking into account that $s_1\to 0$ by Corollary \ref{cors1zero}, $u^-_\e(s_1)=M_1 \to \overline M_1 \in (0, +\infty)$ in view of Proposition \ref{prop6Sect4} and \eqref{bound M1}, $t_1\to \bar t_1 \in (0,1)$, as it follows by combining Proposition \ref{prop6Sect4} and Proposition \ref{tau} (the case $\bar t_1=1$ being excluded by Hopf lemma), then repeating the proofs of \eqref{eq10}, \eqref{eq11}, with $\hat u^-_\e$ replaced by $u^-_\e$, we get that $$\lim_{r \to 0} \bar u^-(r)=\overline M_1, \ \ \lim_{r\to 0}  (\bar u^-)^\prime(r)=0.$$
Hence $\bar u^-$ extends to a $C^1$ radial function near the origin and we easily conclude as in the proof of Proposition \ref{prop6Sect4}.
\end{proof}
\section{Asymptotic analysis of radial sign-changing solutions to \eqref{eq:probM-} with three nodal regions}

To prove Theorem \ref{mainteo2} we could argue by induction starting from $k=2$, nevertheless, for the reader's convenience, we detail the case $k=3$.

Let $u_\e$ be a sign-changing solution of \eqref{eq:probM-} with three nodal regions, let $r_i=r_i(\e)$, $i=1,2$, be the nodal radii, let $s_i=s_i(\e)$, $i=1,2$, be the maximum points of $|u_\e|$ in the second and third nodal region, and denote by $M_i=M_i(\e)$, $i=0,1,2$ the maximum values of $|u_\e|$ in each nodal region.

The following lemma is a trivial consequence of the results obtained in the previous section.
\begin{lemma}\label{lem1Sect5}
As $\e\to 0^+$, we have: $r_1\to 0$, $s_1\to 0$, $\frac{r_1}{r_2}\to0$, $M_0 \to +\infty$ and $$0<\liminf r_2^{\frac{2}{p_\e-1}} M_1\leq\limsup r_2^{\frac{2}{p_\e-1}} M_1<+\infty.$$
\end{lemma}
\begin{proof}
Let us consider the rescaled function
\begin{equation}\label{eqrescaledfunctSect3}
 \tilde u_{\e,2} (x)=r_2^{\frac{2}{p_\e-1}} u_\e(r_2x), \ \ x \in B_{\frac{1}{r_2}}.
\end{equation}
By construction the restriction  of $\tilde u_{\e,2}$ to $B$ is a sign-changing solution to \eqref{eq:probM-} with exactly two nodal regions, and thus we can apply the results of Sect. 4. In particular, denoting by $\tilde r_1 \in (0,1)$ the nodal radius of $\tilde u_{\e,2}$, by $\tilde s_1$ the maximum point of $|\tilde u^-_{\e,2}|$, and setting $\tilde M_0:=\tilde u_{\e,2}(0)$, $\tilde M_1:=|\tilde u_{\e,2}(\tilde s_1)|$, as $\e\to 0^+$ we have:

\begin{itemize}
\item[a)] $\tilde r_1=\frac{r_1}{r_2}\to 0$,
\item[b)] $\tilde s_1=\frac{s_1}{r_2}\to 0$,
\item[c)] $\tilde M_0=r_2^{\frac{2}{p_\e-1}} M_0 \to +\infty$,
\item[d)]  $0<\liminf r_2^{\frac{2}{p_\e-1}} M_1\leq\limsup r_2^{\frac{2}{p_\e-1}} M_1<+\infty$.
\end{itemize}
\end{proof}
 
Next we study the asymptotic behavior of the function $\tilde u_{\e,2}$ defined in \eqref{eqrescaledfunctSect3} in its third nodal region, which is the annulus $A_{1,\frac{1}{r_2}}$. To this end we set $\tilde s_2:=\frac{s_2}{r_2}$, $\tilde M_2:= r_2^{\frac{2}{p_\e-1}} M_2$, which are, respectively, the maximum point and the maximum value of $\tilde u_{\e,2}$ achieved in $A_{1,\frac{1}{r_2}}$.

\begin{proposition}\label{prop1Sect5}
 $\displaystyle\limsup_{\e\to0}\tilde M_2<+\infty$.
\end{proposition}
\begin{proof}
Let us consider the energy-like functionals $H_\lambda:[s_1,t_1] \to \R$, $H_\Lambda:[t_1,s_2]\to \R$ defined by
\begin{equation*}
\begin{split}
H_\lambda(r)&:=\frac{(u_\e^\prime(r))^2}{2}+\frac{|u_\e(r)|^{p_\e+1}}{\lambda(p_\e+1)},\\ H_\Lambda(r)&:=\frac{(u_\e^\prime(r))^2}{2}+\frac{|u_\e(r)|^{p_\e+1}}{\Lambda(p_\e+1)},
\end{split}
\end{equation*}
where $t_1$ is the only point contained in the interval $(r_1,r_2)$ such that $u_\e^{\prime\prime}(t_1)=0$.
Exploiting the ODE in \eqref{eq:initvalprob-}, taking into account that $u_\e^{\prime\prime}\geq0$, $u_\e^{\prime}\geq0$ in $[s_1,t_1]$ and $u_\e^{\prime\prime}\leq0$, $u_\e^{\prime}\geq0$ in $[t_1,s_2]$, we easily check that $H_\lambda$ and $H_\Lambda$ are decreasing. Hence, since $\lambda \leq \Lambda$, we infer that
$$
H_\lambda(s_1)\geq H_\lambda(t_1)\geq H_\Lambda(t_1)\geq H_\Lambda(s_2),
$$
which gives
\begin{equation}\label{eq2-2/10}
\frac{M_1^{p_\e+1}}{\lambda(p_\e+1)}\geq\frac{M_2^{p_\e+1}}{\Lambda(p_\e+1)}.
\end{equation}
From this we get that
$$
r_2^{\frac{2}{p\e-1}}M_1\geq\left(\frac{\lambda}{\Lambda}\right)^{\frac{1}{p_\e+1}}r_2^{\frac{2}{p\e-1}}M_2
$$
and using Lemma \ref{lem1Sect5} we conclude.
\end{proof}

\begin{lemma}\label{prop3Sect5}
$\displaystyle\liminf_{\e\to0}r_2> 0$.
\end{lemma}
\begin{proof}
Assume by contradiction that there exists a sequence $\e \to 0^+$ such that $r_2\to 0$. Consider the restriction to $A_{\frac{r_1}{r_2},\frac{1}{r_2}}$ of the rescaled function $\tilde u_{\e,2}$ defined in \eqref{eqrescaledfunctSect3}. Since $r_2\to 0$ and thanks to Lemma \ref{lem1Sect5} the limit domain of $A_{\frac{r_1}{r_2},\frac{1}{r_2}}$ is $\RN\setminus\{0\}$. Thanks to Lemma \ref{lem1Sect5}, Proposition \ref{prop1Sect5} and elliptic regularity theory we infer that, up to a further subsequence, $\tilde u_{\e,2}\to \tilde u$ in $C^2_{loc}(\RN\setminus\{0\})$, for some radially symmetric function $\tilde u$ satisfying 
$$
-\Mm(D^2 u)=|u|^{p_-^*-1}u \ \ \hbox{in} \ \RN.
$$
Moreover, in view of Theorem \ref{teo1Sect4}, it holds that $\tilde u<0$ in $B$. Hence, since $\tilde u=0$ on $\partial B$, by Hopf's Lemma we get that $\tilde u^\prime(1)>0$. Therefore, for $r>1$ the function $\tilde u=\tilde u(r)$ is a solution (defined and positive in the whole $(1,+\infty)$) to \eqref{eq:initvalprob-} with $p=\psm$, $\alpha=\tilde u^\prime(1)$, but this contradicts \cite[Theorem 1.1]{GIL}. 
\end{proof}

\begin{corollary}\label{cor 3-10}
$\displaystyle\liminf_{\e\to0}\tilde M_2>0$.
\end{corollary}
\begin{proof}
Arguing as in Remark \ref{rem:firsteigen} we have $\tilde M_2^{p_\e-1}\geq \lambda^+_1(-\Mm;A_{1,1/r_2})$ and the conclusion follows from Lemma \ref{prop3Sect5}.
\end{proof}

Finally, summing up, we can describe the asymptotic behavior of $u_\e$.
 
\begin{theorem}\label{teo1Sect5}
Let $u_\e$ be a radial sign-changing solution to \eqref{eq:probM-} with three nodal regions. Then, up to a subsequence, as $\e \to 0^+$ we have that $u_\e \to \bar u$ in $C^2_{loc}(\overline B\setminus\{0\})$, where $\bar u$ is a radial sign-changing solution of \eqref{eq:probM-crit2} having two nodal regions.
\end{theorem}

\begin{proof}
We first observe that, as a consequence of Lemma \ref{lem1Sect5} and Lemma \ref{prop3Sect5}, we have that $M_1$ is uniformly bounded, and bounded away from zero. The same holds for $M_2$ in view of \eqref{eq2-2/10} and Corollary \ref{cor 3-10}. 
 Moreover from Lemma \ref{lem1Sect5} we know that $r_1\to 0$.
Hence, the restriction of $u_\e$ to $A_{r_1,1}$ is uniformly bounded and by standard regularity theory, up to a subsequence, $u_\e \to \bar u$ in $C^2_{loc}(\overline B\setminus\{0\})$, for some radially symmetric function $\bar u$ satisfying \eqref{eqprobwholeball}.

We claim that $\bar u$ is sign-changing with exactly two nodal regions. To prove this we first notice that since $M_2$ is uniformly bounded we have $r_2\not\to 1$, otherwise $\lambda_1(-\Mm; A_{r_2,1})\to+\infty$ and from the inequality $M_2^{p_\e-1} \geq \lambda_1(-\Mm; A_{r_2,1})$ (see Remark \ref{rem:firsteigen}) we would obtain a contradiction. Hence, from this and Lemma \ref{prop3Sect5} we infer that $r_2\to \bar r_2 \in (0,1)$. 

Now, since the restriction to the unit ball of $\tilde u_{\e,2}$ defined in \eqref{eqrescaledfunctSect3} is a sign-changing solution of \eqref{eq:probM-} with two nodal regions, then from Theorem \ref{teo1Sect4} we get that $\tilde u_{\e,2}$ converges in $C_{loc}^2(\overline B\setminus\{0\})$ to the unique negative radial solution $\tilde u$ of \eqref{eq:probM-crit2}. Moreover, since $r_2\to \bar r_2 \in (0,1)$ we deduce that $\bar u(x) = \bar r_2^{-\frac{2}{p_-^*-1}}\tilde u(\bar r_2^{-1} x)$ for $x\in B_{\bar r_2}\setminus\{0\}$. Therefore, $\bar u$ extends to a smooth function near the origin which is a radial solution of \eqref{eq:probM-crit2} and such that $\bar u<0$ in $B_{\bar r_2}$ and $\bar u(\bar r_2)=0$. Hence $\bar u^\prime(\bar r_2)>0$ and thus we easily deduce that $\bar u>0$ in $A_{\bar r_2,1}$. Moreover we have $s_2\to \bar s_2$, for some $\bar s_2$ such that $\bar r_2<\bar s_2 <1$. In fact, $s_2 \to \bar r_2$ cannot happen because $\bar u^\prime(\bar r_2)>0$, while $s_2\not\to 1$ because $\bar u>0$ in $A_{\bar r_2,1}$, $\bar u(1)=0$ and thus $\bar u^\prime(1)<0$.


\end{proof}

\section{Asymptotic analysis of radial sign-changing solutions to \eqref{eq:probM-} with $k$ nodal regions}

In this section we prove Theorem \ref{mainteo2}.

\begin{proof}[Proof of Theorem \ref{mainteo2}]
The case $k=2$ has been proved in Theorem \ref{teo1Sect4}. For $k\geq 3$ we argue by induction on $k$. The case $k=3$ is given by Theorem \ref{teo1Sect5}. Then, assuming the assertion true for a solution with $k$ nodal domains, let $u_\e$ be a radial sign-changing solution of \eqref{eq:probM-} with $(k+1)$ nodal regions. Consider the rescaled function
 \begin{equation}\label{eqrescaledfunctSect5}
 \tilde u_{\e,k} (x)=r_k^{\frac{2}{p_\e-1}} u_\e(r_kx), \ \ x \in B_{\frac{1}{r_k}}.
\end{equation}
Denoting with $\tilde r_i$, $\tilde s_i$, $\tilde M_i$ the corresponding quantities for $\tilde u_{\e,k}$, by construction, we have $\tilde M_i=r_k^{\frac{2}{p_\e-1}}M_i$, $\tilde s_i=\frac{s_i}{r_k}$ for $i=0,\ldots, k$ and $\tilde r_i=\frac{r_i}{r_k}$, for $i=1,\ldots, k$.

Now, since the restriction of $\tilde u_{\e,k}$ to the unit ball is a solution to \eqref{eq:probM-} with $k$ nodal regions then, by the inductive hypothesis, up to a subsequence, as $\e\to 0^+$ we get that $\tilde M_0 \to +\infty$, $\tilde r_1 \to 0$, $\tilde s_1\to 0$ and
$\tilde r_i \to \overline{\tilde r_i}$, $\tilde s_i \to \overline{\tilde s_i}$, for $i=2,\ldots,k-1$, where $0<\overline{\tilde r_2}<\overline{\tilde s_2}<\ldots<\overline{\tilde r_{k-1}}<\ \overline{\tilde s_{k-1}}<1$ and $\tilde M_i \to \overline{\tilde M_i}$, for some  positive numbers $\overline{\tilde M_i}$, for $i=1,\ldots,k-1$. Moreover
 $\tilde u_{\e,k}\to \overline{\tilde u_{k}}$, in $C^2_{loc}(\overline B\setminus\{0\})$, where $\overline{\tilde u_k}$ is a radial sign-changing solution to \eqref{eq:probM-crit2} with $k-1$ nodal regions.
 
  As an immediate consequence we obtain that $M_0 \to +\infty$, $r_1\to 0$, $s_1\to 0$, since $M_0 > \tilde M_0=r_k^{\frac{2}{p_\e-1}} M_0$ and $0<r_1<s_1<\tilde s_1=\frac{s_1}{r_k}$. We divide the remaining part of the proof in three steps:\\
  
\noindent\textbf{Step 1}:  $\tilde M_{k}$ is bounded.\\

Let us consider the energy-like functionals
\begin{equation*}
\begin{split}
\tilde H_\lambda(r)&:=\frac{(\tilde u_{\e,k}^\prime(r))^2}{2}+\frac{|\tilde u_{\e,k}(r)|^{p_\e+1}}{\lambda(p_\e+1)},\\ \tilde H_\Lambda(r)&:=\frac{(\tilde u_{\e,k}^\prime(r))^2}{2}+\frac{|\tilde u_{\e,k}(r)|^{p_\e+1}}{\Lambda(p_\e+1)}.
\end{split}
\end{equation*}
If $k+1$ is even, then $\tilde H_\lambda(r)$ is monotone decreasing in $[\tilde s_{k-1},\tilde s_k]$. Then $\tilde H_\lambda(\tilde s_{k-1})\geq \tilde H_\lambda(\tilde s_{k})$ from which we deduce that 
\begin{equation}\label{eq1 3-10}
\tilde M_{k-1}\geq \tilde M_k.
\end{equation} 
If $k+1$ is odd then $\tilde H_\lambda(r)$ is monotone decreasing in $[\tilde s_{k-1},\tilde t_{k-1}]$, while $\tilde H_\Lambda(r)$ is monotone decreasing in $[\tilde t_{k-1},\tilde s_k]$. Hence, using that $\lambda\leq\Lambda$, we have
$$
\tilde H_\lambda(\tilde s_{k-1})\geq\tilde H_\lambda(\tilde t_{k-1})\geq\tilde H_\Lambda(\tilde t_{k-1})\geq\tilde H_\Lambda(\tilde s_k).
$$
From this we easily deduce that 
\begin{equation}\label{eq2 3-10}
\tilde M_{k-1}\geq \left(\frac\lambda\Lambda\right)^\frac{1}{p_\e+1}\tilde M_k.
\end{equation} 

Exploiting the inductive hypothesis we deduce from \eqref{eq1 3-10}-\eqref{eq2 3-10} that $\tilde M_k$ is bounded. \\

\noindent\textbf{Step 2}: $\displaystyle\liminf_{\e\to0}r_k>0$.\\

Assume by contradiction that $r_k \to 0$, for some subsequence $\e\to 0$. Then, the limit domain of $A_{\tilde r_1, \frac{1}{r_k}}$ is $\R^N\setminus\{0\}$ and, in view of Step 1 and the inductive hypothesis, we have that the restriction $\tilde u_{\e,k}\big|_{A_{\tilde r_1, \frac{1}{r_k}}}$ is uniformly bounded. Hence $\tilde u_{\e,k} \to \tilde u_k$ in $C^2_{loc}(\R^N\setminus\{0\})$, where $\tilde u_k$ is a radial solution of
\begin{equation}\label{eq:probM-critRN}
-\Mm(D^2 u)=|u|^{p^*_- -1}u \ \ \ \hbox{in} \ \R^N\setminus\{0\}.
\end{equation}
Moreover, $\tilde u_k$ is non-trivial and sign-changing because by construction we have $\tilde u_k\equiv \overline{\tilde u_k}$ in $\overline B\setminus\{0\}$, where $\overline{\tilde u_k}$ is the limit of $\tilde u_{\e,k}\big|_B$\,. In particular, $\tilde u_k^\prime(1)\neq 0$, $\tilde u_k(1)=0$ and thus it cannot happen that $\tilde u_k\equiv 0$ in $\R^N \setminus \overline B$. Therefore, the restriction of $\tilde u_k$ to $\R^N \setminus \overline B$ is a constant-sign radial solution of
 \begin{equation}\label{eq:probextdomSect6}
\begin{cases}
-\Mm(D^2 u)=|u|^{p_-^*-1}u & \hbox{in} \ \R^N \setminus \overline B,\\
\qquad  \ \ \ \ \ u=0 & \hbox{on} \ \partial B.
\end{cases}
\end{equation}
This contradicts \cite[Theorem 1.1]{GIL}.\\

\noindent\textbf{Step 3}: conclusion.\\

In view of Step 1 and Step 2, since $\tilde M_k=r_k^{\frac{2}{p_\e-1}} M_k$, we infer that $M_k$ is bounded. Moreover it is bounded away from zero. Indeed, arguing as in Remark \ref{rem:firsteigen} we get that $M_k^{p_\e-1}\geq\lambda_1^+\left(-\Mp;B\right)$ if $k+1$ is odd, or $M_k^{p_\e-1}\geq\lambda_1^+\left(-\Mm;B\right)$ if $k+1$ is even.

 Up to a subsequence we then have $M_k \to \bar M_k$,  $r_k \to \bar r_k$, for some $\bar M_k\in (0,+\infty)$, $\bar r_k \in (0,1]$. Arguing as in the proof of Theorem \ref{teo1Sect5}, taking the restriction of $u_\e$ to $A_{r_k, 1}$, exploiting that $M_k$ is bounded and that $\lambda^+_1(-\Mpm; A_{r_k,1}) \to +\infty$ if $r_k \to 1$, we infer that $\bar r_k \neq 1$. 

Now, using the definitions of $\tilde M_i$, $\tilde r_i$, $\tilde s_i$ for $i=2,\ldots, k-1$ and the results proved in the first part of the proof, we conclude that $M_i\to \bar M_i=\bar r_k^{- \frac{2}{\psm-1}} \overline{\tilde M_i}\in (0,+\infty)$, $r_i\to \bar r_i={\overline{ \tilde r_i}}{\bar r_k}$, $s_i\to \bar s_i= {\overline{ \tilde s_i}}{\bar r_k}$ for $i=2,\ldots,k-1$, and it holds $$0<\bar r_2<\bar s_2<\ldots<\bar r_{k-1} < \bar s_{k-1}.$$

Summing up, since the restriction of $u_\e$ to $A_{r_1, 1}$ is uniformly bounded and $\bar r_k \neq 1$ we deduce that, up to a subsequence, $u_\e \to \bar u$ in $C^2_{loc}(\overline B\setminus\{0\})$, where $\bar u$ is a non-trivial radial sign-changing solution of
\begin{equation}\label{eq:probM-crit2pd}
\begin{cases}
-\Mm(D^2 u)=|u|^{p^*_- -1}u & \hbox{in} \ B\setminus\{0\},\\
\qquad  \ \ \ \ \ u=0 & \hbox{on} \ \partial B,
\end{cases}
\end{equation}
with $k$ nodal regions. Since $\bar r_k \in (0,1)$, ${\overline{ \tilde s_{k-1}}} \in (0,1)$ we infer that $\bar s_{k-1}<\bar r_k$. Moreover, denoting by $\bar s_k$ the limit point of $s_k$, from the regularity of $\bar u$ in compact subsets of $\overline B\setminus\{0\}$ and since $\bar M_k\neq 0$, we deduce that it cannot happen that $\bar s_k=\bar r_k$. 
Moreover, since $u_\e'(s_k)=0$ and $|\bar u'(1)|>0$ in view of the Hopf's Lemma, we infer that $s_k \to 1$ cannot happen.
Therefore the nodes and the extrema of $\bar u$ are ordered in the following way
$$0<\bar r_2<\bar s_2<\ldots<\bar r_{k-1} < \bar s_{k-1}<\bar r_k<\bar s_k<1,$$
as expected.

At the end,  arguing as in the proof of Theorem \ref{teo1Sect5} we see that $\bar u(x)=\bar r_k^{- \frac{2}{\psm-1}}\overline{\tilde u_k}(\bar r_k x)$, for $x \in B_{\bar r_k}\setminus\{0\}$. Hence $\bar u$ extends to a $C^2$ radial function near the origin and it is a sign-changing solution of \eqref{eq:probM-crit2} having $k$ nodal regions. This completes the proof of the inductive step.
\end{proof}

\section{Asymptotic analysis of radial sign-changing solutions to \eqref{eq:probM+} with two nodal regions}


 Let $u_n$ be a sequence of radial sign-changing solutions to \eqref{eq:probM+} with two nodal domains, where $p_n=\psc-\e_n$, $\e_n\searrow 0$ as $n\to +\infty$. The first result is about the behavior of $M_n:=\|u_n\|_\infty$ as $n\to +\infty$. We set $M_0=M_0(n):=\|u_n^+\|_{\infty}=u_n(0)$ and $M_1=M_1(n):=\|u_n^-\|_{\infty}$, and we denote by $r_1=r_1(n)$ the nodal radius and by $s_1=s_1(n)$ the minimum point.

\begin{proposition}\label{prop1Sect7}
We have $M_n \to +\infty$, as $n\to +\infty$.
\end{proposition}
\begin{proof}
Arguing as in Remark \ref{rem:firsteigen} we infer that
\beq\label{eq:limitationm0m1}
\begin{split}
 M_0^{\psce-1} &\geq  \lambda^+_1(-\Mp;B_{r_1})\geq \lambda^+_1(-\Mp;B),\\
 M_1^{\psce-1} &\geq  \lambda^+_1(-\Mm;A_{r_1,1})\geq  \lambda^+_1(-\Mm;B),
\end{split}
\eeq
which readily implies that $M_0$ and $M_1$ are uniformly bounded from below away from zero. Since $M_n = \max\left\{M_0,M_1\right\}$ the same holds for $M_n$.

Now, assume that $M_n\to\overline M\in(0,+\infty)$ for some subsequence. Then, from elliptic regularity estimates we deduce that $u_n$ is uniformly bounded in $C^2(\overline B)$. In addition, in view of \eqref{eq:limitationm0m1} it cannot happen that $r_1\to0$ or $r_1\to1$, otherwise $\lambda^+_1(-\Mp,B_{r_1}) \to +\infty$ or $ \lambda^+_1(-\Mm;A_{r_1,1})\to +\infty$ and by \eqref{eq:limitationm0m1} this would imply that one between $M_0$ and $M_1$ blows-up, contradicting the uniform boundedness of $u_n$. 

Therefore $r_1\not\to 0$, $r_1\not\to1$ and by regularity estimates, up to a further subsequence, as $n\to +\infty$, we have $u_n \to \bar u$  in $C^2(\overline B)$, where $\bar u$ is a radial sign-changing solution of \eqref{probllimupodd} and $\bar u(0)>0$ (because of \eqref{eq:limitationm0m1}), but this contradicts Proposition \ref{prop:nonexistscsolMP}.
\end{proof}

\begin{proposition}\label{prop2Sect7}
We have:
\begin{itemize}
\item[i)]  $\displaystyle0<\liminf_{n\to+\infty}r_1 M_0^{\frac{\psce-1}{2}}\leq\limsup_{n\to+\infty}r_1 M_0^{\frac{\psce-1}{2}}<+\infty$;
\item [ii)] $\displaystyle0<\liminf_{n\to+\infty}\frac{M_0}{M_1}\leq\limsup_{n\to+\infty}\frac{M_0}{M_1}<+\infty$;
\item[iii)] $M_0 \to +\infty$, $M_1\to +\infty$;
\item[iv)] $r_1 \to 0$;
\item[v)] the rescaled function 
\beq\label{globalrescaledfunct2}
\tilde u_n(x):= r_1^{\frac{2}{\psce-1}} u_n\left(r_1x\right), \ x \in B_{1/r_1},
\eeq
converges, up to a subsequence, in $C^2_{loc}(\RN)$ to a non-trivial radial sign-changing solution of
\beq\label{limitproblemrescaledfunct}
-\Mp(D^2 u)=|u|^{\psc-1}u \ \ \hbox{in} \ \RN.
\eeq
\item [vi)] $s_1\to 0$.
\end{itemize}

\end{proposition}
\begin{proof}
If along a subsequence it holds that $r_1 M_0^{\frac{\psce-1}{2}}\to +\infty$, then, defining the rescaled function $$\hat u_n(x):=\frac{1}{M_0} u_n\left(\frac{x}{M_0^{\frac{\psce-1}{2}}}\right), \ x\in B_{M_0^{\frac{\psce-1}{2}}},$$
we have $\hat u_n(0)=1$ and by the usual argument we infer that $\hat u_n^+ \to \hat u$ in $C^2_{loc}(\RN)$, where $\hat u$ is a positive radial solution of 
\begin{equation}\label{eq:probM+poslimite}
-\Mp(D^2 u)=u^{\psc} \ \ \hbox{in} \ \RN.
\end{equation}
On the other hand, since $\psc<\psp$ then by \cite[Theorem 1.1]{FQ} we get that \eqref{eq:probM+poslimite} does not have  positive radial solutions, and this gives a contradiction.

If along a subsequence it holds that $r_1 M_0^{\frac{\psce-1}{2}} \to 0$, arguing as in Remark \ref{rem:firsteigen} for $\hat u_n$, we would have $$1 \geq \lambda^+_1(-\Mp,B_{r_1 M_0^{\frac{\psce-1}{2}}})\to +\infty$$ which is a contradiction. This completes the proof of i).\\

To prove ii) we show that $\frac{M_0}{M_1} \to 0$ and $\frac{M_0}{M_1} \to +\infty$ cannot occur along any subsequence.

Let us consider the functionals $$H_\lambda(r):=\frac{(u_n^\prime)^2(r)}{2} + \frac{u_n^{p_n+1}}{\lambda(\psce+1)},\ \ r \in [0,t_0],\ \ \ \ H_\Lambda(r):= \frac{(u_n^\prime)^2(r)}{2} + \frac{|u_n|^{\psce+1}}{\Lambda(\psce+1)},\ \ r \in [t_0,s_1],$$
 where $t_0 \in (0,r_1)$ is the only radius such that $u_n^{\prime\prime}<0$ in $(0,t_0)$ and  $u_n^{\prime\prime}>0$ in $(t_0,r_1)$. Exploiting the ODE satisfied by $u_n$ in $[0,s_1]$, we check that $H_\lambda$ is decreasing in $[0,t_0]$ and $H_\Lambda$ is decreasing in $[t_0,s_1]$. Moreover, since $\lambda \leq \Lambda$ and $u_n(t_0)>0$ we infer that $H_\lambda(t_0)\geq H_\Lambda(t_0)$. Summing up we get that
$$H_\lambda(0)\geq H_\lambda(t_0)\geq H_\Lambda(t_0)\geq H_\Lambda(s_1), $$ 
and since $H_\lambda(0)=\frac{M_0^{\psce+1}}{\lambda(p_n+1)}$, $H_\Lambda(s_1)=\frac{M_1^{\psce+1}}{\Lambda(p_n+1)}$ we deduce that
\beq\label{eqM1M2}
\frac{M_0}{M_1}\geq \left(\frac{\lambda}{\Lambda}\right)^{\frac{1}{\psce+1}}.
\eeq
From \eqref{eqM1M2} it follows that 
$$
M_n\leq\left(\frac{\Lambda}{\lambda}\right)^{\frac{1}{\psce+1}}M_0
$$
and this implies 
\begin{equation}\label{tt}
\lim_{n\to+\infty}M_0=+\infty.
\end{equation}

Assume now that, for some subsequence, $\frac{M_0}{M_1} \to +\infty$ and consider again the rescaled function $\hat u_n$.
By construction $\hat u_n(0)=1$, $\left\|\hat u_n\right\|_\infty\leq 1$ and as before, up to a subsequence as $n\to +\infty$, we have $\hat u_n \to \hat u$ in $C^2_{loc}(\RN)$, where $\hat u$ is a non-trivial radial solution to \eqref{limitproblemrescaledfunct} satisfying $\hat u(0)=1$.

From i) there exists $c_1>0$ such that $\hat u>0$ in $B_{c_1}$ and $\hat u=0$ on $\partial B_{c_1}$. Moreover $\hat u^\prime(c_1)<0$. On the other hand, taking into account that  $\frac{M_0}{M_1} \to +\infty$, we deduce that $\hat u\equiv 0$ in $\R^N\backslash\overline B_{c_1}$, since for any fixed $x$ such that $|x| > c_1$ and for all sufficiently large $n$ we have
$$|\hat u_n(x)|=\left|  \frac{1}{M_0} u_n\left(\frac{x}{M_0^{\frac{\psce-1}{2}}}\right)\right| \leq \frac{M_1}{M_0}.$$
Passing to the limit as $n\to +\infty$ we obtain $\hat u(x)=0$, which contradicts the $C^1$ regularity of $\hat u$.  

Statement iii) is an immediate consequence of \eqref{tt} and ii), while iv) directly follows from i) and iii). \\

Let us prove v) and vi).  From i) and ii) the rescaled function $\tilde u_n$ in \eqref{globalrescaledfunct2} is uniformly bounded, and by iii) the limit of $B_{1/r_1}$ is $\RN$. Moreover, by construction $\tilde u_n\big|_B= v_{\psce,+}$ is the only radial positive solution of \eqref{eq:probM+}. Hence, up to a subsequence as $n\to +\infty$, we have $\tilde u_n \to \tilde u$ in $C^2_{loc}(\RN)$, where $\tilde u$ is a non-trivial radial solution of \eqref{limitproblemrescaledfunct}, and this proves v). In particular $\tilde u>0$ in $B$, $\tilde u=0$ on $\partial B$ and the function $\tilde u^-=\tilde u^-(r)$ coincides, for $r>1$, with the unique maximal solution of \eqref{eq:initvalprob-} with $p=\psc$, $\alpha= -(\tilde u^+)^\prime(1)$.
Now arguing as in the proof of the first part of Proposition \ref{prop sigma}, we infer that $\frac{s_1}{r_1}$ is bounded. Since $r_1\to0$ by iv), we conclude that $s_1\to0$ as well. 

\end{proof}
 
In the next result we state and prove a uniform upper bound that will be crucial in the sequel (see Sect. 10).
\begin{proposition}\label{prop:unifupperbound}
Let $\tilde u_n^-$ be the negative part of the rescaled function $\tilde u_n$ defined in \eqref{globalrescaledfunct2}.
 Then, up to a subsequence as $n\to +\infty$, $\tilde u_n^- \to W_-$ in $C^2_{loc}(\RN\setminus B)$, where $W_-$  is the only positive radial fast decaying solution of \eqref{eq:domextmmpsc}. Moreover, there exist two positive constants $C$, $K$ (independent on $n$) such that for all sufficiently large $n$ it holds
 \beq\label{unifuppbound}
\tilde u_n^-(r)\leq  \frac{C}{\left(r^2-(\tilde t_1)^2 + K\right)^{\frac{\Nm-2}{2}}} \qquad \hbox{for} \ r \in \left[\tilde t_1,\frac{1}{r_1}\right],
\eeq
where $\tilde t_1 \to \bar t_1\in (1,+\infty)$ and  $\bar t_1$ is the only radius where $W_-=W_-(r)$ changes concavity.
\end{proposition}

\begin{proof}
In view of Proposition \ref{prop2Sect7}-v), up to a subsequence as $n\to +\infty$, we have in particular that $\tilde u_n^- \to \tilde u^-$ in $C^2_{loc}(\RN\setminus B)$, where $\tilde u^-$ is a positive radial  solution of \eqref{eq:domextmmpsc}. We claim that $\tilde u^-$ is fast decaying. By (ii) of Theorem \ref{fastMm} this is equivalent to show that $\alpha_-^*(\psc)=(\tilde u^-)^\prime(1)$.  To prove this, we observe that $\tilde u_n^-=\tilde u_n^-(r)$ is a solution of \eqref{eq:initvalprob-} with $p=\psce$, $\alpha=\alpha(n)=(\tilde u^-_n)^\prime(1)=-v'_{\psce,+}(1)$ and such that $\tilde u^-_n(1/r_1)=0$, and thus, by definition of  $\alpha^*_-$, we have
$$ \alpha^*_-(\psce) \leq (\tilde u^-_n)^\prime(1).$$
Then, passing to the limit as $n\to +\infty$, exploiting Proposition \ref{lem:asalphastar} and taking into account that $\tilde u_n \to \tilde u$ in $C^2_{loc}(\RN)$, we deduce that
$$\alpha_-^*(\psc)\leq (\tilde u^-)^\prime(1).$$
On the other hand, since $\tilde u^-=\tilde u^-(r)$ is defined and positive in the whole $(1, +\infty)$, then, by (i) of Theorem \ref{fastMm}, we infer that $\alpha_-^*(\psc)\geq(\tilde u^-)^\prime(1)$. Hence, $\tilde u^-=W_-$ is the only positive radial fast decaying solution of \eqref{eq:domextmmpsc}.

For \eqref{unifuppbound}, we notice that it is exactly inequality \eqref{eq:uniformupperboundalpha}  obtained  in the proof of Proposition \ref{lem:asalphastar}. Indeed, denoting by $\tilde t_1$ the only radius where $\tilde u_n^-=\tilde u_n^-(r)$ changes concavity and setting $v_n(r):=\tilde u_n^-(r)$, we have that $v_n$ satisfies \eqref{ODEunifuppbound} in $\left[\tilde t_1, \frac{1}{r_1}\right]$, with $p_n=\psc - \e_n$. In the present case,  $p_n\to p^{**}_+$, which still satisfies  $\psc>\psm>\frac{\Nm+2}{\Nm-2}$. Then the proof follows verbatim, taking into account that $v_n \to W_-$ in $C^2_{loc}([1,+\infty)$, and that $\tilde t_1\to \bar t_1$, $v_n(\tilde t_1) \to W_-(\bar t_1)$, as $n\to +\infty$, where $\bar t_1$ is the only radius where $W_-=W_-(r)$ changes concavity.
\end{proof}

We conclude this section by studying the limiting behavior of $u_n$.
\begin{proposition}\label{PropConvSolM+}
Up to a subsequence, as $n\to +\infty$, it holds that $u_n \to 0$ in $C^2_{loc}(\overline B\setminus\{0\})$.
\end{proposition}
\begin{proof}
The proof is carried out along the same line of \cite[Theorem 1.1, ii)]{BGLP}.  By elliptic estimates, it is sufficient to show that for any fixed $\rho\in(0,1)$  $$\left\|u_n\right\|_{\infty,A_{\rho}}\to0\qquad\text{as $n\to +\infty$,}$$   where $A_\rho=\left\{x\in\overline B:\,|x|\geq\rho\right\}$ and $\|\cdot\|_{\infty, A_\rho}$ denotes the $L^\infty$-norm in $A_\rho$. \\
Using \eqref{unifuppbound} and the definition \eqref{globalrescaledfunct2} of $\tilde u_n$ we have 
$$r_1^{\frac{2}{\psce-1}} u_n^-(r_1r)\leq C \frac{1}{\left(r^2-\frac{t_1^2}{r_1^2} + K\right)^{\frac{\Nm-2}{2}}} \qquad \hbox{for} \ r \in \left[\frac{t_1}{r_1},\frac{1}{r_1}\right),$$
or equivalently
$$r_1^{\frac{2}{\psce-1}} u_n^-(r)\leq C \frac{1}{\left(\frac{r^2}{r_1^2}-\frac{t_1^2}{r_1^2} + K\right)^{\frac{\Nm-2}{2}}}, \qquad \hbox{for} \ r \in \left[t_1,1\right)$$
for some positive constants $C$ and $K$ independent on $n$.\\
Now, since $\frac{t_1}{r_1}\to\bar t_1\in(1,+\infty)$ by Proposition \ref{prop:unifupperbound} and $r_1\to0$ in view of Proposition \ref{prop2Sect7}-iv), we infer that $t_1\to0$. Hence for sufficiently large $n$ we have $t_1^2\leq\frac{\rho^2}{2}$  and 
\begin{equation*}
\left\|u_n\right\|_{\infty,A_{\rho}}=\left\|u^-_n\right\|_{\infty, A_{\rho}}\leq\frac{C}{\left(\frac{\rho^2}{2}\right)^{\frac{\Nm-2}{2}}}r_1^{\Nm-2-\frac{2}{\psce-1}}.
\end{equation*}
 Then the conclusion follows  since $r_1\to0$ and $\Nm-2-\frac{2}{\psce-1}\to\Nm-2-\frac{2}{p^{**}_+-1}>0$.

\end{proof}

\begin{remark}\label{rem:newphenomenon}
By Proposition \ref{prop2Sect7} and Proposition \ref{PropConvSolM+} we have that the maximum and the minimum of $u_n$, namely $M_0$ and $-M_1$, blow up at the same rate and the minimum point $s_1$ converges to the maximum point which is zero. Thus we have concentration of the positive and negative part at the same point.

This is a new phenomenon, as compared with the classical Lane-Emden problem in which case, whenever the rate of blow-up of the positive and negative part is the same, the two nodal regions separate and the concentration points of the negative and positive part are different (see \cite{BEP}).
\end{remark}

\section{Asymptotic analysis of radial sign-changing solutions to \eqref{eq:probM+} with three nodal regions}\label{M+3}

Having proved in the previous section Theorem \ref{mainteo4} when $k=2$ we could argue by induction to get the general result. However, since passing from even to odd the statement changes, we prefer to detail the proof for $k=3$, for the reader convenience.

Let $u_n$ be a sign-changing solution of \eqref{eq:probM+} with three nodal regions, let $r_i=r_i(n)$, $i=1,2$, be the nodal radii, let $s_i=s_i(n)$, $i=1,2$, be the maximum points of $|u_n|$ in the second and third nodal region, and denote by $M_i=M_i(n)$, $i=0,1,2$ the maximum values of $|u_n|$ in each nodal region.

\begin{proposition}
Up to a subsequence, as $n\to +\infty$, we have:
 $M_0\to +\infty$,  $M_1\to +\infty$, $r_i\to 0$, $s_i\to 0$ for $i=1,2$, $M_{2} \to \bar M$, for some $\bar M>0$, and $u_n \to \bar v$ in $C^2_{loc}(\overline{B}\setminus\{0\})$, where $\bar v$ is the unique positive radial solution of \eqref{probllimupodd}.
\end{proposition}
\begin{proof}
Let us consider the rescaled function
\begin{equation}
\tilde u_n(x):=r_2^{\frac{2}{\psce-1}}u_n(r_2x),\quad x\in B_{1/r_2}.
\end{equation}
Then the restriction of $\tilde u_n$ to the unit ball $B$ is the radial sign-changing solution of \eqref{eq:probM+} with two nodal regions so that the results of Sect. 7 apply. In particular, up to a subsequence, as $n\to+\infty$, we have:
\begin{equation}\label{1010eq1}
\frac{r_1}{r_2}\to0\;, \qquad r_2^{\frac{2}{\psce-1}}M_i\to+\infty\quad\text{for $i=0,1$}
\end{equation}
and
\begin{equation}\label{1510eq1}
\frac{M_0}{M_1}=\frac{r_2^{\frac{2}{\psce-1}}M_0}{r_2^{\frac{2}{\psce-1}}M_1}\to c_0,
\end{equation}
for some positive constant $c_0$.
From  \eqref{1010eq1} we  deduce that 
\begin{equation*}
\begin{split}
M_i&\to+\infty\quad\text{for $i=0,1$}\\
r_1&\to0.
\end{split}
\end{equation*}
Moreover
$$
 \tilde H(r)=\frac{(\tilde u_{n}^\prime(r))^2}{2}+\frac{|\tilde u_{n}(r)|^{\psce+1}}{\lambda(\psce+1)}
$$
is monotone decreasing in $[1,\frac{s_2}{r_2}]$, hence 
\begin{equation}
\frac{{\left(r_2^{\frac{2}{\psce-1}}M_2\right)}^{\psce+1}}{\lambda(\psce+1)}\leq\frac{(\tilde u_{n}^\prime(1))^2}{2}.
\end{equation}
Since $\tilde u_{n}^\prime(1)\to0$, then
\begin{equation}\label{1015eq2}
{r_2^{\frac{2}{\psce-1}}M_2}\to0
\end{equation}
 and, using the lower bound $M_2^{\psce-1} \geq\lambda^+_1(-\Mp;B)$ (see \eqref{eq:limitationm0m1}), we also deduce that $r_2\to0$ and, as a consequence, $s_1\to0$. Putting together \eqref{1010eq1} for $i=1$ and \eqref{1015eq2} we have
\begin{equation}
\frac{M_1}{M_2}=\frac{r_2^{\frac{2}{\psce-1}}M_1}{r_2^{\frac{2}{\psce-1}}M_2}\to+\infty.
\end{equation}
Arguing as in Proposition \ref{prop sigma} and Corollary \ref{cors1zero} we also deduce that $s_2\to0$.
If $M_2$ is bounded from above, then $u_n\big|_{A_{r_2,1}} \to \bar u$ in $C^2_{loc}(\overline B\setminus\{0\})$ for some radial positive function $\bar u$. Since $s_2\to 0$, then as in the proof of Theorem \ref{teo1Sect4} we obtain that $\bar u$ extends to the unique  positive radial solution of \eqref{probllimupodd}, as we aim to prove. Hence to complete the proof it remains to show that $\limsup_{n\to +\infty}M_2<+\infty$. On the contrary, let us assume that, along some sequence $n\to +\infty$, $M_2\to+\infty$  and consider the rescaled function
$$\hat u_{n}(x):=\frac{1}{M_2} u_n\left(\frac{x}{M_2^{\frac{\psce-1}{2}}}\right), \ \ x\in A_{r_2 M_2^{\frac{\psce-1}{2}}, M_2^{\frac{\psce-1}{2}}}.$$ 
Since the limit domain of $\hat u_n$  is $\R^N\backslash\left\{0\right\}$ by  \eqref{1015eq2}, we can   argue exactly as in Proposition \ref{prop6Sect4} to deduce that $\hat u_{n} \to \hat u$ in $C^2_{loc}(\RN\setminus\{0\})$, where $\hat u$ can be extended to a non-trivial radial positive solution to \eqref{eq:probM+poslimite}. This is clearly a contradiction since $\psc<\psp$. Hence the only possibility is that $M_2\to \bar M$ for some positive constant $\bar M$ as we wanted to show.
\end{proof}

\section{Asymptotic analysis of radial sign-changing solutions to \eqref{eq:probM+} with $k$ nodal regions }

In this section we prove Theorem \ref{mainteo4}. 
\begin{proof}[Proof of Theorem \ref{mainteo4}.]
We argue by induction on $k$. The steps $k=2,3$ have been proved respectively in Sect. 7 and Sect. 8. So let us assume that the assertion holds true for solutions with $k$ nodal regions and let $u_n$ be a radial sign-changing solution to \eqref{eq:probM+} with $k+1$ nodal regions.\\

\noindent\textbf{Case 1:} if $k+1$ is even, then the restriction to $B$ of the rescaled function
 \begin{equation*}
 \tilde u_{n,k} (x):=r_k^{\frac{2}{\psce-1}} u_n(r_kx), \ \ x \in B_{1/{r_k}},
\end{equation*}
is a solution to \eqref{eq:probM+} having $k$ nodal regions, with $k$ odd. Hence, exploiting the inductive hypothesis and the definition of $\tilde u_{n,k}$, we infer that, up to a subsequence as $n\to +\infty$,  $M_0 r_k^{\frac{2}{\psce-1}} \to +\infty,\ldots,M_{k-2} r_k^{\frac{2}{\psce-1}} \to +\infty$, and $M_{k-1} r_k^{\frac{2}{\psce-1}} \to \bar M$, for some $\bar M>0$. This readily implies that $M_0\to +\infty, \ldots,M_{k-2}\to +\infty$. Moreover, by inductive hypothesis, we also have $\frac{r_i}{r_k} \to 0$, $\frac{s_i}{r_k} \to 0$ for $i=1,\ldots,k-1$, which easily implies that ${r_i}\to 0$, ${s_i}\to 0$, for $i=1,\ldots,k-1$.

Finally, exploiting again the inductive hypothesis we deduce that
$\frac{M_{2j}}{M_{2j+1}}=\frac{M_{2j}r_k^{\frac{\psce-1}{2}}}{M_{2j}r_k^{\frac{\psce-1}{2}}} \to c_j$ for $j=0,\ldots,\frac{k-3}{2}$, with $c_j$ positive constants, and that
$\frac{M_{2j+1}}{M_{2j+2}}=\frac{M_{2j+1}r_k^{\frac{\psce-1}{2}}}{M_{2j+2}r_k^{\frac{\psce-1}{2}}} \to +\infty$ for $j=0,\ldots,\frac{k-3}{2}$.

Now, arguing as in the proof of \eqref{eqM1M2} we have that
\beq\label{eqrationextremalvalue}
\frac{M_{k-1}}{M_{k}}\geq \left(\frac{\lambda}{\Lambda}\right)^{\frac{1}{\psce+1}}.
\eeq
Let us show that $$\limsup_{n\to+\infty}\frac{M_{k-1}}{M_{k}}<+\infty.
$$
Assume by contradiction that, for a subsequence, $\frac{M_{k-1}}{M_{k}}\to +\infty$ and consider the restriction of $\tilde u_{n,k}$ to the annulus $A_{{r_{k-1}}/{r_k}, 1/r_k}$. Since $M_{k-1} r_k^{\frac{2}{\psce-1}} \to \bar M$, for some $\bar M>0$ then from  \eqref{eqrationextremalvalue} we deduce that $\tilde u_{n,k}\big|_{A_{{r_{k-1}}/{r_k}, 1/r_k}}$ is uniformly bounded. We claim that $r_k\not\to 1$. Otherwise, since $\tilde u_{\e,k}\big|_{A_{1, 1/r_k}}$ is  uniformly bounded and $r_k \to 1$ we would have $\lambda_1(-\Mm; A_{1, 1/r_k}) \to +\infty$, and from \eqref{bound M1} we would have a contradiction. Hence $r_k \to \bar r \in  (0,1)$ or $r_k \to 0$, and recalling that by inductive hypothesis $\frac{r_{k-1}}{r_k} \to 0$, we infer that $\tilde u_{n,k}\big|_{A_{r_{k-1}/r_k, 1/r_k}} \to \tilde u$ in $C^2_{loc}(\Pi)$, where $\tilde u$ is non-trivial  (because by inductive hypothesis it coincides with $\bar u$ in $B$) and the limit domain $\Pi$ is either $\RN\setminus\{0\}$ or $\overline{B}_{\frac{1}{\bar r}}\setminus\{0\}$. Since we are assuming that $\frac{M_{k-1}}{M_{k}}\to +\infty$, arguing as in the proof of Proposition \ref{prop2Sect7}, ii), we deduce that $\tilde u\equiv 0$ in $\Pi \cap \{x\in \RN; \ |x|>1\}$ contradicting the regularity of $\tilde u$. Hence the only possibility is 
\beq\label{eq:ratiomkm1mk}
\frac{M_{k-1}}{M_{k}}\to c_{(k-1)/2},
\eeq
for some $c_{(k-1)/2}>0$. Let us also show that $r_k \to 0$. Indeed, since $
\tilde u_{n,k}\big|_{A_{{r_{k-1}}/{r_k}, 1/r_k}} \to \tilde u$ in $C^2_{loc}(\Pi)$, where $\tilde u$ can be extended to a non-trivial sign-changing solution to $-\Mp u=|u|^{\psc-1}u$ in $B_{\bar r}$, with homogenous Dirichlet boundary condition if $\Pi=\overline{B}_{\bar r}\setminus\{0\}$, or $-\Mp u=|u|^{\psc-1}u$ in $\RN$ if $\Pi=\RN\setminus\{0\}$. Now, in view of Proposition \ref{prop:nonexistscsolMP} the first case cannot occur and and thus we infer that $r_k \to 0$. In particular, since $M_{k-1}r_k^{\frac{2}{\psce-1}} \to \bar M$, it follows that $M_{k-1} \to +\infty$ and from \eqref{eq:ratiomkm1mk} we infer that $M_k \to +\infty$.

Moreover, since 
$M_{k}s_k^{\frac{2}{\psce-1}}$ is bounded from above (by arguing as in the first part of Proposition \ref{prop sigma}) and $M_k \to +\infty$, then $s_k \to 0$.

To conclude the proof it remains to show that $u_n \to 0$ in $C^2_{loc}(\overline{B}\setminus\{0\})$. The proof of this fact is identical to that of Proposition \ref{PropConvSolM+} with minor modifications. In particular, taking into account that $r_k,s_k\to 0$ and $u_n$ is negative in the last nodal component, then, for any fixed $\rho \in (0,1)$ and for all sufficiently large $n$ we obtain 
$$\left\|u_n\right\|_{\infty, A_{\rho}}=\left\|u^-_n\right\|_{\infty, A_{\rho}}\leq\frac{C}{\left(\frac{\rho^2}{2}\right)^{\frac{\Nm-2}{2}}}r_k^{\Nm-2-\frac{2}{\psce-1}},
$$
for some positive constant $C=C(N,\lambda,\Lambda)$ independent on $n$.  This completes the proof when $k+1$ is even.\\

\noindent\textbf{Case 2:} if $k+1$ is odd we consider the restriction to $B$ of the rescaled function
 \begin{equation}\label{eqrescaledfunctSect8}
 \tilde u_{n,k} (x):=r_k^{\frac{2}{\psce-1}} u_n(r_kx), \ \ x \in B_{1/{r_k}},
\end{equation}
which is a solution to \eqref{eq:probM+} having $k$ nodal regions, with $k$ even. Hence, from the inductive hypothesis, up to a subsequence, as $n\to +\infty$, we get that $M_0 r_k^{\frac{2}{\psce-1}} \to +\infty,\ldots,M_{k-1} r_k^{\frac{2}{\psce-1}} \to +\infty$ and thus we infer that $M_0\to +\infty, \ldots,M_{k-1}\to +\infty$. Exploiting again the inductive hypothesis, we have $\frac{r_i}{r_k} \to 0$, $\frac{s_i}{r_k} \to 0$ for $i=1,\ldots,k-1$, which implies that ${r_i}\to 0$, ${s_i}\to 0$, for $i=1,\ldots,k-1$. Moreover $\frac{M_{2j}}{M_{2j+1}}=\frac{M_{2j} r_k^{\frac{\psce-1}{2}}}{M_{2j+1}r_k^{\frac{\psce-1}{2}}} \to c_j$ for $j=0,\ldots,\frac{k-2}{2}$ and $c_j$ positive constants, $\frac{M_{2j+1}}{M_{2j+2}}=\frac{M_{2j+1} r_k^{\frac{\psce-1}{2}}}{M_{2j+2}r_k^{\frac{\psce-1}{2}}} \to+\infty$, for $j=0,\ldots,\frac{k-4}{2}$.\\

Repeating exactly the same arguments of Sect. \ref{M+3}, in the case of three nodal regions, we infer that $r_k,s_k\to0$, $\frac{M_{k-1}}{M_k}\to+\infty$
and that $M_k\to\bar M$, where $\bar M$ is a positive constant. From this we also deduce that $u_n$ converges in $C^2_{loc}(\overline B\backslash\left\{0\right\})$ to the unique positive solution of \eqref{probllimupodd}. The proof is complete.
\end{proof}

\section{Energy of solutions}

Let $\Omega$ be a bounded radial domain in $\R^N$, i.e. $\Omega$ is either a ball or an annulus centered at the origin. Then the radial coordinate $r$ will belong either to $[0,R)$, $R>0$, if $\Omega$ is the ball $B_R$, or to the interval $(a,b)$, $0<a<b$, if $\Omega$ is the annulus $A_{a,b}$.

We consider the space of radial functions in $\Omega$ which have constant sign and change convexity only once. More precisely we define
\begin{eqnarray*}
X_\Omega&:=&\left\{u \in C_{rad}^2(\overline\Omega); \ |u|>0 \ \text{and}\ \exists\ \varrho=\varrho(u) \in (a,b)\ [\text{resp.} \ \varrho \in (0,R) \ \text{if} \ \Omega \ \text{is a ball}]  \right.\\
&& \ \ \text{such that} \ u^{\prime\prime}(\varrho)=0, \ u^{\prime\prime}(r) <0 \ \text{for} \ r\in(a,\varrho) \ \text{and} \ u^{\prime\prime}(r)>0 \ \text{for} \ r\in(\varrho,b), \\
&&  \ \ \left. \text{or} \ u^{\prime\prime}(r) >0 \ \text{for} \ r\in(a,\varrho) \ \text{and} \ u^{\prime\prime}(r)<0 \ \text{for} \ r\in(\varrho,b), \right.\\
&&\ \  \left. [\text{resp.} \  r\in(0,\varrho)  \ \text{and} \ r\in(\varrho,R) \ \text{if} \ \Omega \ \text{is a ball}] \right\}.
\end{eqnarray*}
Next, for an exponent $p>1$ and for any function $u \in X_\Omega$ we consider the radial weight:

\begin{equation}\label{radialweight}
g_{u,p}(x):=\begin{cases}
[\varrho(u)]^{\gamma(p)} & \text{if} \ |x|=r\leq \varrho(u),\\
|x|^{\gamma(p)} & \text{if} \ |x|>\varrho(u),
\end{cases}
\end{equation}
with $\gamma(p):=2\left(\frac{p+1}{p-1}\right)-N$, and define the weighted energy
\begin{equation}\label{weightedenergy}
E_{p,\Omega}(u):= \int_\Omega |u(x)|^{p+1} g_{u,p}(x) \ dx.
\end{equation}
It is elementary to check that $E_{p,\Omega}$ is invariant under the  scaling $u_\alpha(x)=\alpha u(\alpha^{\frac{p-1}{2}}x)$ (see \cite[Proof of Thorem 1.2]{BGLP}), i.e. 
\begin{equation}\label{energyinvariance}
E_{p,\Omega}(u)=E_{p,\Omega_\alpha}(u_\alpha)
\end{equation}
with  $\Omega_\alpha=\alpha^{-\frac{p-1}{2}}\Omega$.

We observe that if $u$ is a solution of \eqref{eq:probMgen} with $k$ nodal regions, then the restrictions $u^m$ to each nodal region $\Omega^m$, for $m=1,\ldots,k$, belong to the space $X_{\Omega^m}$.
Therefore we can consider the energy of each function $u^m$ in the corresponding nodal region $\Omega^m$:
\begin{equation}\label{restrweightedenergy}
E_{p,\Omega^m}(u^m):= \int_{\Omega^m} |u^m(x)|^{p+1} g_{u^m, p}(x) \ dx, \ \ \text{for} \ m=1,\ldots,k,
\end{equation}
and define the total energy of the solution $u$ in the ball $B$ as
\begin{equation}\label{wtotalenergy}
E_{p}^T(u):=\sum_{m=1}^k E_{p, \Omega^m}(u^m).
\end{equation}
A similar energy can be defined for any positive (fast decaying) radial solution $U_{\pm}$ of the critical equation in $\R^N$
\begin{equation}\label{criticalproblemRN}
-\Mpm(D^2u)=u^{\pspm} \ \ \text{in} \ \R^N.
\end{equation}
We denote it by $E^*(U_{\pm})$, i.e.
\begin{equation}\label{weightedenergyRN}
E^*(U_{\pm}):= \int_{\R^N} |U_\pm(x)|^{\pspm +1} g_{U_\pm}^{*}(x) \ dx,
\end{equation}
where
\begin{equation}\label{radialweightRN}
g_{U_\pm}^*(x):=\begin{cases}
[\varrho(U_\pm)]^{\gamma^*} & \text{if} \ |x|=r\leq \varrho(U_\pm),\\
|x|^{\gamma^*} & \text{if} \ |x|=r>\varrho(U_\pm),
\end{cases}
\end{equation}
with $\gamma^*:=2\left(\frac{\pspm+1}{\pspm-1}\right)-N$.

Note that, by the invariance of the energy with respect the usual scaling, $\Sigma^*_\pm:=E^*(U_{\pm})$ is a constant depending  only on $\lambda,\Lambda, N$. \\

We now prove Theorem \ref{mainteo3} (we refer to Theorem \ref{mainteo2} for the notations).
\begin{proof}[Proof of Theorem \ref{mainteo3}]
Let us first consider the restriction of the solution $u_\e$ to the first nodal region $\Omega^1_\e=B_{r_1}$. We denote it by $u_\e^1$. The function
$$
\tilde u_\e^1(x):=r_1^{\frac{2}{p_\e-1}}u_\e^1(r_1x),\quad x\in B,
$$
is the unique positive  solution of 
\begin{equation*}
\begin{cases}
-\Mm(D^2 u)=u^{p_\e} & \hbox{in} \ B,\\
\qquad  \ \ \ \ \ u=0 & \hbox{on} \ \partial B.
\end{cases}
\end{equation*}
By the scaling invariance \eqref{energyinvariance}, we have 
$$
E_{p_\e, \Omega^1_\e}(u_\e^1)= E_{p_\e, B}(\tilde u_\e^1).
$$
By \cite[Theorem 1.2]{BGLP} we immediately deduce that, as $\e\to0$,
$$
E_{p_\e, \Omega^1_\e}(u_\e^1)\to\int_{\R^N} (U_-)^{p^*_-+1} g^*_{U_-} \ dx=\Sigma^*_-.
$$
This gives the first contribution to the limit of the total energy in \eqref{eq:limenergy}.

On the other hand we recall that the nodal radii $r_1,r_2,\ldots, r_{k-1}$ converge respectively to $0$, $\bar r_2$,\ldots, $\bar r_{k-1}$, where $\bar r_2$,\ldots, $\bar r_{k-1}$ are the nodal radii of the limit function $\bar u$ given by Theorem \ref{mainteo2}. 

Thus, by the convergence of $u_\e \to \bar u$ in $C^2_{loc}(\overline B\setminus\{0\})$ we have that the restriction $u_\e^m$ of $u_\e$ to its nodal region $\Omega_\e^m$, $m=2,\ldots,k$, converges to the restriction of $\bar u$ to the corresponding nodal region, i.e.:
\begin{equation*}
u_\e^2 \to\bar u_1,\ \ldots,\ u_\e^m \to \bar u_{m-1}.
\end{equation*}
Then, using also that $M_i \to \bar M_i$, we have
\begin{eqnarray*}
E_{p_\e, \Omega_\e^2}(u_\e^2)&=&\int_{\Omega_\e^2} |u_\e^2|^{p_\e+1} g_{u_\e^2} \ dx \to \int_{B_{\bar r_2}} |\bar u_1|^{p^*_- +1} g^*_{\bar u_1} \ dx\\
&\ldots&\\
E_{p_\e, \Omega_\e^m}(u_\e^m)&=&\int_{\Omega_\e^m} |u_\e^m|^{p_\e+1} g_{u_\e^m} \ dx \to \int_{A_{\bar r_{m-1}, \bar r_{m}}} |\bar u_m|^{p^*_- +1} g^*_{\bar u_{m-1}} \ dx,
\end{eqnarray*}
where $\bar r_m=1$ if $m=k$.
Thus the assertion \eqref{eq:limenergy} holds and the proof is complete.
\end{proof}

We now study the limit energy of a family of sign-changing solution to \eqref{eq:probM+} having $k$ nodal regions, as $n\to +\infty$ and prove Theorem \ref{mainteo5}. We denote by $E^{**}(W_-)$ the energy of the only radial positive fast decaying solution to \eqref{eq:domextmmpsc}, namely
\beq\label{eq:energyW-}
E^{**}(W_-)=\int_{\RN\setminus B}  |W_-(x)|^{\psc+1} g^{**}_{W_-}(x) \ dx, 
\eeq
where
\begin{equation}\label{radialweightWRN}
g_{W_-}^{**}(x):=\begin{cases}
[\varrho(W_-)]^{\gamma^{**}} & \text{if} \ 1\leq |x|=r\leq \varrho(W_-),\\
|x|^{\gamma^{**}} & \text{if} \ \ \ \ \ \ \ |x|=r>\varrho(W_-),
\end{cases}
\end{equation}
with $\gamma^{**}:=2\left(\frac{\psc+1}{\psc-1}\right)-N$. Since $W_-$ is fast decaying easily we have:

\begin{lemma}\label{lemmaEfinita}
The energy $E^{**}(W_-)$ is finite.
\end{lemma}
\begin{proof}
Since $W_-$ is fast decaying we can find $C>0$ and $\bar r>1$ such that 
\beq\label{eq:fastdecayW}
W_-(r)\leq \frac{C}{r^{\Nm-2}} \ \ \text{for } \ r>\bar r.
\eeq
Up to choosing a larger $\bar r$ we can assume without loss of generality that $\bar r>\varrho(W_-)$ and thus
\begin{eqnarray*}
E^{**}(W_-) &=&\int_{\RN\setminus B}  |W_-(x)|^{\psc+1} g^{**}_{W_-}(x) \ dx\\
&=&\int_{1}^{\varrho(W_-)} |W_-(r)|^{\psc+1} [\varrho(W_-)]^{2\left(\frac{\psc+1}{\psc-1}\right)-N} r^{N-1} \ dr\\
&& + \int_{\varrho(W_-)}^{+\infty} |W_-|^{\psc+1} r^{2\left(\frac{\psc+1}{\psc-1}\right)-N} r^{N-1} \ dr\\
&=& (I) + (II).
\end{eqnarray*} 
Clearly $(I)$ is finite. For $(II)$, exploiting \eqref{eq:fastdecayW} we have
\begin{eqnarray*}
(II) &\leq& \int_{\varrho(W_-)}^{\bar r} |W_-(r)|^{\psc+1} r^{2\left(\frac{\psc+1}{\psc-1}\right)-N} r^{N-1} \ dr + C  \int_{\bar r}^{+\infty} r^{-(\Nm-2)(\psc+1) + 2\left(\frac{\psc+1}{\psc-1}\right)-1} \ dr\\
&=& (III) + (IV).
\end{eqnarray*}
Now, $(III)$ is finite and so is $(IV)$ by a straightforward computation because $\psc>\psm>\frac{\Nm}{\Nm-2}$.
\end{proof}
We now prove Theorem \ref{mainteo5} (for the notations we refer to Theorem \ref{mainteo4} and the beginning of this section).

\begin{proof}[Proof of Theorem \ref{mainteo5}]
We argue by induction on $k\geq 2$. We begin with the basic step $k=2$.

Assume that $u_n$ is a sign-changing solution to \eqref{eq:probM+} having two nodal regions and let $r_1$ be the node of $u_n$. Then $\Omega^1_n=B_{r_1}$ and $\Omega^2_n=A_{r_1,1}$. We consider the rescaled function $\tilde u_n$ defined in \eqref{globalrescaledfunct2}. By construction $\tilde u_n^+\big|_B$ coincides with the unique positive solution of \eqref{eq:probM+}, i.e. $\tilde u_n^+\big|_B=v_{\psce,+}$ and it is uniformly bounded (see Proposition \ref{prop2Sect7}). Hence $\tilde u_n^+\big|_B \to \bar v$ in $C^2(\overline B)$, where $\bar v$ is the unique positive  solution of  \eqref{probllimupodd}. Hence, exploiting the scaling invariance of the energy and passing to the limit as $n\to +\infty$ we get that
\beq\label{eqlimitenergy1}
E_{\psce, \Omega^1_n}(u^1_n)=E_{\psce, B}(\tilde u^+_n)=\int_B  |\tilde u_n^+|^{\psce+1} g_{\tilde u_n^+} \ dx \to \int_B  |\bar v|^{\psc+1} g_{\bar v} \ dx=E_{\psc, B}(\bar v).
\eeq
For $u^2_n$ we have $\Omega^2_n=A_{r_1,1}$ and exploiting the scaling invariance and the definition of the energy we have
\begin{eqnarray*}
E_{\psce, \Omega^2_n}(u^2_n)=E_{\psce, A_{1,1/r_1}}(\tilde u^-_n)&=&\int_{A_{1,\tilde t_1}} |\tilde u_n^-|^{\psce+1} [\tilde t_1]^{\gamma(\psce)} \ dx + \int_{A_{\tilde t_1,1/r_1}} |\tilde u_n^-|^{\psce+1} |x|^{\gamma(\psce)} \ dx,\\
&=& (I) + (II),
\end{eqnarray*}
where $\tilde t_1=\varrho(\tilde u^-_n)$. From Proposition \ref{prop:unifupperbound} we have $\tilde u_n^-\big|_{A_{1,1/r_1}} \to W_-$ in $C^2_{loc}(\RN\setminus B)$, where $W_-$ is the unique radial positive fast decaying solution of \eqref{eq:domextmmpsc}. We claim that we can pass to limit under the integral sign in $(I)$ and $(II)$. Indeed, for  $(I)$ it is obvious because $\tilde u_n$ is uniformly bounded and $\tilde t_1\to \bar t_1 \in (0,+\infty)$, where $\bar t_1=\varrho(W_-)$, while for $(II)$,  taking into account \eqref{unifuppbound} and that $\psce>\psm>\frac{\Nm}{\Nm-2}$, we easily obtain
\begin{eqnarray*}
 |\tilde u_n^-(r)|^{\psce+1} r^{2\frac{\psce+1}{\psce-1}-N} r^{N-1} \leq C r^{-(\Nm-2)(\psce+1)+2\frac{\psce+1}{\psce-1}-1} \leq C r^{-1-\delta} \ \ \text{for $r\geq \tilde t_1$,}
\end{eqnarray*}
for some $C>0$, $\delta>0$ independent on $n$.
Hence, by Lebesgue's dominated convergence theorem we can pass to the limit under the integral sign in (II) and thus we conclude that
\beq\label{eqlimitenergy2}
\begin{array}{lll}
\displaystyle \lim_{n\to +\infty} E_{\psce, \Omega^2_n}(u^2_n)&=&\displaystyle \int_{A_{1,\bar t_1}} |W_-(x)|^{\psc+1} [\bar t_1]^{\gamma(\psc)} \ dx+ \int_{\R^N\backslash B_{\bar t_1}} |W_-(x)|^{\psc+1} |x|^{\gamma(\psc)} \ dx\\[12pt]
&=& E^{**}(W_-)=\Sigma^{**}_+.
\end{array}
\eeq
Combining \eqref{eqlimitenergy1} and \eqref{eqlimitenergy2} we complete the proof of the basic step.\\

Now let us prove the inductive step. Let $u_n$ be a radial solution to \eqref{eq:probM+} with $k+1$ nodal regions and consider the rescaled function $\tilde u_{k,n}$ defined by \eqref{eqrescaledfunctSect8}. Clearly, by invariance under this scaling and by definition we easily have
\beq\label{totalenergym+}
E_{\psce}^T(u_n)= E_{\psce}^T(\tilde u_{k,n}\big|_B)+ E_{\psce, A_{1,1/r_k}}(\tilde u_{k,n}\big|_{A_{1,1/r_k}}).
\eeq
Since the restriction $\tilde u_{k,n}\big|_B$ is a radial sign-changing solution of \eqref{eq:probM+} with $k$ nodal regions then by the induction hypothesis we infer that
\beq\label{totenergyeven0}
\lim_{n\to +\infty}E_{\psce}^T(\tilde u_{k,n}\big|_B)=\begin{cases}
  \frac{k}{2} E_{\psc, B}(\bar v) +  \frac{k}{2}\Sigma^{**}_+& \hbox{if $k$ is even},\\[4pt] 
 \frac{k+1}{2} E_{\psc, B}(\bar v) +  \frac{k-1}{2}\Sigma^{**}_+  & \hbox{if $k$ is odd}.
\end{cases}
\eeq
For the second term of \eqref{totalenergym+}, if $k$ is even then $u_n$ has an odd number of nodal components and by Theorem \ref{mainteo4} we have $r_k\to 0$, $u_n^{k+1}=u_n\big|_{A_{r_k,1}}$ is uniformly bounded and  $u_n^{k+1} \to \bar v$ in $C^2_{loc}(\overline B\setminus\{0\})$. Therefore, by the usual invariance under scaling and exploiting these properties we get that
\beq\label{totenergyeven}
 \lim_{n\to +\infty} E_{\psce, A_{1,1/r_k}}(\tilde u_{k,n}\big|_{A_{1,1/r_k}})= \lim_{n\to +\infty} E_{\psce,A_{r_k,1}}(u_n^{k+1})=E_{\psc, B}(\bar v).
\eeq
If $k$ is odd, then applying Theorem \ref{mainteo4} to $\tilde u_{k,n}\big|_B$ we infer that $\tilde u_{k,n}\big|_B \to \bar v$ in $C^2_{loc}(\overline B\setminus\{0\})$. Hence, arguing as in the proof of Proposition \ref{prop:unifupperbound} we have $\tilde u_n^-\big|_{A_{1,1/r_k}} \to W_-$ in $C^2_{loc}(\RN\setminus B)$ and $\tilde t_k \to \bar t_k \in (1, +\infty)$, $\bar t_k=\varrho(W_-)$. Then, as in the proof of the basic step $k=2$ and exploiting an analogous uniform upper bound (the proof is the same as that of Proposition \ref{prop:unifupperbound} with minor changes), we get that
\beq\label{totenergyodd}
\lim_{n\to +\infty} E_{\psce, A_{1,1/r_k}}(\tilde u_{k,n}\big|_{A_{1,1/r_k}})=\Sigma^{**}_+.
\eeq
At the end, combining \eqref{totalenergym+}--\eqref{totenergyodd} we obtain \eqref{eq:mainteo5}.

\end{proof}

\end{document}